\documentclass[final,leqno,onefignum,onetabnum]{siamltexmm}

 \usepackage{bm}
\usepackage[numbers]{natbib}
\usepackage{url}
\usepackage[intlimits]{amsmath}
\usepackage{graphicx}
\usepackage{fancyhdr}
\usepackage{amsfonts}
\usepackage{amssymb}
\usepackage{array}
\newcolumntype{P}[1]{>{\centering\arraybackslash}p{#1}}

\usepackage{float}
\usepackage[caption=false]{subfig}
\usepackage{subfiles}
\usepackage{microtype}

\usepackage{hyperref}
\usepackage[noabbrev]{cleveref}

\crefname{figure}{Figure}{Figures}
\Crefname{figure}{Figure}{Figures}
\begin{document}

\title{An adaptive preconditioner for steady incompressible flows
%\title{Pinning and Depinning in the Swift--Hohenberg equation with time-periodic forcing
%\thanks{ }
}
%\title{Pinning and depinning in the Swift--Hohenberg equation with time-periodic forcing}
\author{
C\'edric Beaume 
\thanks{Department of Applied Mathematics, University of Leeds, Leeds LS2 9JT, UK (\email{ced.beaume@gmail.com}).}
}
%\date{\today}
\maketitle
\newcommand{\slugmaster}{%
\slugger{siads}{xxxx}{xx}{x}{x--x}}%slugger should be set to juq, siads, sifin, or siims

\begin{abstract}
This paper describes an adaptive preconditioner for numerical continuation of incompressible Navier--Stokes flows. The preconditioner maps the identity (no preconditioner) to the Stokes preconditioner (preconditioning by Laplacian) through a continuous parameter and is built on a first order Euler time-discretization scheme. The preconditioner is tested onto two fluid configurations: three-dimensional doubly diffusive convection and a reduced model of shear flows. In the former case, Stokes preconditioning works but a mixed preconditioner is preferred. In the latter case, the system of equation is split and solved simultaneously using two different preconditioners, one of which is parameter dependent. Due to the nature of these applications, this preconditioner is expected to help a wide range of studies.
\end{abstract}



\pagestyle{myheadings}
\thispagestyle{plain}
\markboth{C. Beaume}{A adaptive preconditioner for continuation of fluid flows}

%%%%%%%%%%%%%%%%%%%%%%%%%%%%%%%%%%%%%%%%%%%%%%%%%%%%%%%%%%%%%%%%%%%%%%%
%%%%%%%%%%%%%%%%%%%%%%%%%%%%%%%%%%%%%%%%%%%%%%%%%%%%%%%%%%%%%%%%%%%%%%%
\section{Introduction}
%%%%%%%%%%%%%%%%%%%%%%%%%%%%%%%%%%%%%%%%%%%%%%%%%%%%%%%%%%%%%%%%%%%%%%%
%%%%%%%%%%%%%%%%%%%%%%%%%%%%%%%%%%%%%%%%%%%%%%%%%%%%%%%%%%%%%%%%%%%%%%%

The development of specialized numerical methods and the increase in available computing resources have helped making substantial progress in understanding many nonlinear problems as dynamical systems.
The most basic tool available to that end is time integration which simulates the temporal evolution of an initial condition, thereby emulating an experimental or natural realization.
Time integration provides access to the preferred transient and end state, however, it does not (necessarily) provide access to information regarding the origin of these end states.
One way to understand how these states are formed is to compute unstable solutions.
These states cannot be obtained, or in some rare cases very hardly, using time integration but help provide a complete picture of the dynamical system: these can gain stability or lead to the creation of new solutions or new transients under parametrical changes.
%Another of its limitations is that it is not designed to chart parameter space, as each time a parameter is changed, the simulation has to be reinitialized and the transient simulated.
Numerical continuation has been developed to complement time-integration in that respect and has become an essential part of the toolkit of the nonlinear dynamicist.

Pionnered by Keller \cite{Keller77}, these methods compute steady solutions of a system of ordinary differential equations (ODE) and their evolution as the values of the parameter of the problem are changed.
Their aim is to {\it continue} a fixed point in parameter space in order to draw its {\it branch} and unfold the {\it bifurcation diagram} explaining its formation.
Continuation methods consist in a two-step algorithm comprising a prediction phase using previous iterates along the branch and a correction phase involving a fixed point method \cite{Allgower80,Allgower85,Allgower03,Seydel09}.
Due to their nature, these methods are capable of computing exact solutions regardless of their stability and provide information on the effect of parametric changes on a solution.
Numerical continuation became a very popular tool, broadly used in many different fields \cite{Seydel87,Henry00,Krauskopf07} and a myriad of packages have been developped and released in the open domain \cite{Engelborghs01,Kuznetsov03,Clewley07,Doedel08,Uecker14}.

The area of fluid dynamics has seen much progress with the help of continuation methods.
Intricate pattern formation problems have been elucidated such as that of Rayleigh--B\'enard convection rolls in cartesian \cite{Torres13,Torres14}, cylindrical \cite{Boronska10} and spherical shell geometries \cite{Feudel11}.
More complicated physics has been tackled, such as doubly diffusive convection \cite{Bergeon02} and free surface binary fluid convection \cite{Bergeon98} with similar success.
Spatially localized pattern formation, involving large aspect-ratio domains, has also been investigated: a collection of spatially localized convective states have been found in two-dimensional large aspect-ratio binary fluid convection \cite{Mercader09,Mercader11}, rotating convection \cite{Beaume13rc1,Beaume13rc2} and magnetoconvection \cite{Lojacono11}.
Despite the successful and reliable use of continuation methods in two-dimensional and small three-dimensional domains, the extention to more complex geometries constitutes a major challenge.
The most noticeable attempts concern doubly diffusive convection in a three-dimensional domain of square cross section and large transverse direction \cite{Beaume13ddd1,Beaume13ddd2} and porous medium convection in domains extended in two directions \cite{Lojacono13}, each of these problems involving $O(10^6)$ degrees of freedom.
These studies involved unreasonably long simulation campaigns and require a certain level of experience in the use of numerical continuation to provide results in an efficient way, if at all.

Another area of fluid dynamics that has benefited from the developments of numerical continuation is that of transition to turbulence.
Shear flows such as plane Couette flow or pipe flow are subcritical flows, i.e., the trivial laminar solution is stable and coexists with turbulence, a state in which the flow displays spatial and temporal complexity, above a threshold value of the parameters.
Meticulous studies of the unstable exact solutions living between both states have provided crucial understanding of transition.
The pioneering discovery of unstable nonlinear solutions in plane Couette flow \cite{Nagata90} drew a lot of attention and started a new research area.
The simplest of these solutions has been thoroughly studied \cite{Wang07,Kawahara12} and new solutions discovered \cite{Clever97,Waleffe03,Schneider08,Gibson09,Melnikov14} that all together provide a comprehensive picture of transitional phenomena.
Similar studies took place in other shear flows and hinted at a common mechanism for transition to turbulence in subcritical shear flows \cite{Faisst03,Wedin04,Pringle07}.
Lastly, recent tours de force involving numerical continuation on large domains revealed families of spatially localized states of different kinds \cite{Schneider10,Chantry14,Gibson14}.

The aim of this paper is to describe an efficient and adaptive way to precondition the Navier--Stokes equation governing steady incompressible flows.
Some basic principles of numerical continuation are summarized in Section \ref{secnumcont}.
In Section \ref{precon}, I describe the preconditioning method, followed in Section \ref{secex} by two examples: three-dimensional coupled convection and shear flows.
A short conclusion terminates the paper.

%%%%%%%%%%%%%%%%%%%%%%%%%%%%%%%%%%%%%%%%%%%%%%%%%%%%%%%%%%%%%%%%%%%%%%%
%%%%%%%%%%%%%%%%%%%%%%%%%%%%%%%%%%%%%%%%%%%%%%%%%%%%%%%%%%%%%%%%%%%%%%%
\section{Numerical continuation}\label{secnumcont}
%%%%%%%%%%%%%%%%%%%%%%%%%%%%%%%%%%%%%%%%%%%%%%%%%%%%%%%%%%%%%%%%%%%%%%%
%%%%%%%%%%%%%%%%%%%%%%%%%%%%%%%%%%%%%%%%%%%%%%%%%%%%%%%%%%%%%%%%%%%%%%%

We consider the simple dynamical system:
\begin{equation}
\label{fu0}
\partial_t {\bf u} = {\bf F}({\bf u},\lambda),
\end{equation}
where $t$ represents time, ${\bf u} \in \mathcal{R}^n$ is the solution vector of dimension $n$, ${\bf F}: \mathcal{R}^{n+1} \rightarrow \mathcal{R}^{n}$ is a nonlinear operator and $\lambda$ the continuation or free parameter.
We seek solutions that satisfy $\partial_t {\bf u} = 0$ or equivalently ${\bf F}({\bf u},\lambda) = 0$.
Note that in case there is more than one parameter, all the parameters but one ($\lambda$) are kept fixed and are included in the operator ${\bf F}$.
A continuation method consists in parameterizing and continuing the branch of solutions of equation (\ref{fu0}): $\mathcal{B}(s) = ({\bf u}(s),\lambda(s))$, where $s$ is the arclength along the branch.
This is done in two steps: a prediction step based on extrapolation of previous results along the branch is created and then converged with fixed point method based on equation (\ref{fu0}).
In the following, I detail simple, programmer-friendly, choices for the prediction step and then describe two different ways of converging them: fixed parameter and pseudo-arclength continuation.
The former is the simplest to implement and is the one used in the exemples in Section \ref{secex}.
It provides an easy guide to get started with continuation.

\subsection{Prediction}

The first step of a continuation method is prediction.
The simplest way to predict consists in a polynomial extrapolation along the branch at a parametric distance $\triangle \lambda$.
If only one solution along the branch $({\bf u}_1,\lambda_1)$ is known, it is used as a predictor $({\bf u}_{2}^p,\lambda_{2}^p)$ of a second solution $({\bf u}_2,\lambda_2)$:
\begin{eqnarray}
&{\bf u}_{2}^p &= {\bf u}_1,\\
&\lambda_{2}^p &= \lambda_1 + \triangle \lambda.
\end{eqnarray}
Upon successful computation of the second solution, a linear extrapolation is used to provide an initial condition $({\bf u}_{3}^p,\lambda_{3}^p)$ for the third point $({\bf u}_3,\lambda_3)$:
\begin{eqnarray}
&{\bf u}_{3}^p &= {\bf u}_2 + \frac{{\bf u}_2 - {\bf u}_1}{\lambda_2 - \lambda_1} \triangle \lambda,\\
&\lambda_{3}^p &= \lambda_2 + \triangle \lambda.
\end{eqnarray}
From the moment three or more solutions are known, there is typically little to gain in increasing the polynomial degree of the approximation and quadratic extrapolation remains standard:
\begin{eqnarray}
\label{arcpredic}
&{\bf u}_{i}^p &= {\bf u}_{i-1} + \frac{{\bf u}_{i-1} - {\bf u}_{i-2}}{\lambda_{i-1} - \lambda_{i-2}} \triangle \lambda + \left( \frac{{\bf u}_{i-1} - {\bf u}_{i-2}}{\lambda_{i-1} - \lambda_{i-2}} - \frac{{\bf u}_{i-2} - {\bf u}_{i-3}}{\lambda_{i-2} - \lambda_{i-3}} \right) \frac{\triangle \lambda + \lambda_{i-1} - \lambda_{i-2}}{\lambda_{i-1} - \lambda_{i-3}} \triangle \lambda,\\
\label{straightpredic}
&\lambda_{i}^p &= \lambda_{i-1} + \triangle \lambda,
\end{eqnarray}
for the prediction $({\bf u}_{i}^p,\lambda_{i}^p)$ for the $i$-th solution $({\bf u}_i,\lambda_i)$.

The method above, although simple to derive and program, fails in the presence of a saddle-node as $d\lambda/ds$ changes sign.
This can easily be fixed by monitoring the approach of a saddle-node and changing the sign of $\triangle \lambda$ when the corrector fails at the approach of a saddle-node.
It is also usual to replace equation (\ref{straightpredic}) by a condition on the arclength $\triangle s$, such as:
\begin{equation}
\label{arcpredic2}
({\bf u}_{i}^p - {\bf u}_{i-1})^2 + (\triangle \lambda)^2 = (\triangle s)^2,
\end{equation}
to predict at a distance $\triangle s$ from solution $(u_i,\lambda_i)$.
Other methods, more sophisticated or tailored to the specific needs or preferences of the user can be generated but are out of the scope of the present article.

\subsection{Correction}

To converge the prediction, one condition has to be added so that equation (\ref{fu0}) yields a well-posed problem.
This first possibility consists in fixing the value of the parameter $\lambda$ and changing the continuation mode at the approach of a saddle-node.
It is called {\it fixed parameter continuation} is explained in the following section.
Another very popular continuation method consists in adding a condition ensuring that the correction is done along a vector orthogonal to the prediction vector.
This method is known under the name of {\it pseudo arc-length continuation} and is described in the subsequent section.

\subsubsection{Fixed parameter continuation}
\label{paramcont}

We use a Newton--Raphson (hereafter Newton) method to converge the prediction obtained in the previous section to the solution of system (\ref{fu0}).
Fixed parameter continuation consists in freezing one quantity, typically $\lambda$, and correct the other ones.
We then set $\lambda_i = \lambda_i^p$ and solve the resulting $n \times n$ system:
\begin{equation}
\label{newton1}
D_{{\bf u}}{\bf F}({\bf u}^j)\delta {\bf u}^j = {\bf F}({\bf u}^j),
\end{equation}
where $D_{{\bf u}}{\bf F}({\bf u}^j)$ is the Fr\'echet derivative of ${\bf F}$ with respect to ${\bf u}$ and the superscript $j$ denotes quantities evaluated at the $j$-th Newton iteration.
The computation of $\delta {\bf u}^j$ is then followed by the correction ${\bf u}^{j+1} = {\bf u}^j - \delta {\bf u}^j$.
Provided the prediction $(u_i^p,\lambda_i^p)$ is close enough to a solution, the Newton method (\ref{newton1}) converges quadratically.

Issues arise at saddle-nodes, as no solution might exist at $\lambda_i^p$.
To anticipate this, an alternate continuation mode is used in which $\lambda$ is free to vary while one element of ${\bf u}$ is kept constant, resulting in an $n \times n$ system again.
If $u_{\{k\}}$ is part of the solution: $\partial_t u_{\{k\}} = F_{\{k\}} ({\bf u}) = 0$, where $u_{\{k\}}$ and $F_{\{k\}}$ denote respectively the $k$-th element of ${\bf u}$ and ${\bf F}$.
We introduce $\hat{\bf u} = {\bf u} \setminus \{u_{\{k\}}\}$ the vector constituted of all elements of ${\bf u}$ but $u_{\{k\}}$ and $\hat{\bf F}({\bf u}) = {\bf F}({\bf u}) \setminus \{F_{\{k\}}({\bf u})\}$ similarly.
Under this rearrangement, continuation of saddle-nodes can be achieved by solving the following system at the $j$-th Newton iteration:
\begin{equation}
\begin{pmatrix}
D_{\hat{\bf u}} \hat{\bf F} (\hat{\bf u}^j,\lambda^j) & D_{\lambda} \hat{\bf F} (\hat{\bf u}^j,\lambda^j) \\
D_{\hat{\bf u}} F_{\{k\}} (\hat{\bf u}^j,\lambda^j) & D_{\lambda} F_{\{k\}} (\hat{\bf u}^j,\lambda^j)
\end{pmatrix}
\begin{pmatrix}
\delta \hat{\bf u}^j \\
\delta \lambda^j
\end{pmatrix}
=
\begin{pmatrix}
\hat{\bf F} ({\bf u}^j,\lambda^j) \\
0
\end{pmatrix}
,
\end{equation}
where the right-hand-side of the last equation is $F_{\{k\}}({\bf u}^j,\lambda^j) = 0$.
The correction reads $\hat{\bf u}^{j+1} = \hat{\bf u}^j - \delta \hat{\bf u}^j$ and $\lambda^{j+1} = \lambda^j - \delta \lambda^j$.

\begin{figure}
\begin{center}
\includegraphics[width=0.8\textwidth]{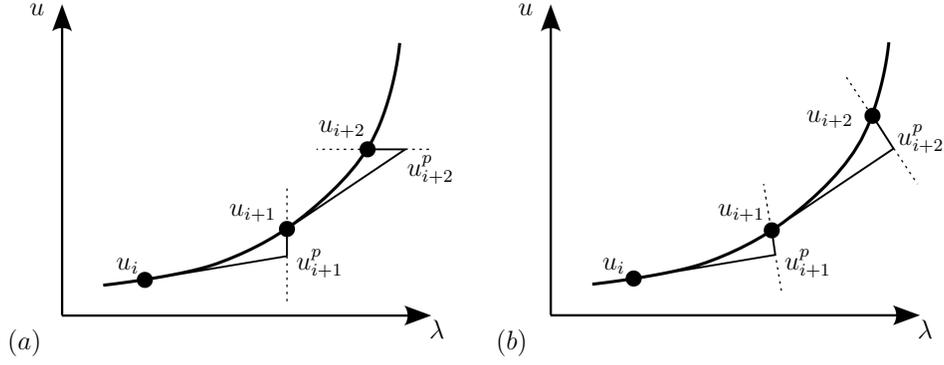}
\end{center}
\caption{Two different types of continuation methods for a scalar unknown $u$ against parameter $\lambda$: fixed parameter continuation (a) and pseudo-arclength continuation (b). The former consists in correcting either at a fixed parameter value (see correction from $u_{i+1}^p$) or at a fixed value of the solution (or one of its components, see correction from $u_{i+1}^p$). The latter consists in correcting on a direction orthogonal to the prediction direction.}
\label{continuation}
\end{figure}

To pass saddle-nodes, a criterion has to be set up to determine when to switch from fixed $\lambda$ to fixed $u_{\{k\}}$ continuation.
This criterion can involve the slope of the branch, for instance:
\begin{equation}
\label{criterion}
\frac{\partial \mathcal{N}({\bf u})}{\partial \lambda} < c,
\end{equation}
where $\mathcal{N}({\bf u})$ represents a norm of ${\bf u}$ and $c$ an arbitrary real constant whose optimal value is problem-dependent.
When criterion (\ref{criterion}) is true, the slope of the branch is gentle and the continuation is done for a fixed $\lambda$.
On the contrary, when criterion (\ref{criterion}) is false, the slope of the branch is significant and might indicate the presence of a saddle-node. 
In that case, switching to a fixed $u_{\{k\}}$ continuation provides better results.
Figure \ref{continuation}(a) provides a sketch of this type of continuation for a scalar solution $u$, highlighting both the fixed parameter correction from $u_{i+1}^p$ and the fixed $u$ correction for $u_{i+2}^p$.
%%%%%%%%%%%%%%%%%%%%%%%%%%
%%%%%%%%%%%%%%%%%%%%%%%%%%

The prediction/correction loop just described is usually optimized through variable prediction distances.
When the prediction distance (either $\triangle \lambda$ or $\triangle s$) is too small, the correction is very simple.
It takes one or two Newton iteration to converge to the desired accuracy and many small trivial steps are made along the branch, wasting time in the process.
When the prediction distance is too large, the correction is tedious and can take many difficult Newton iterations, jump over to another branch or simply fail.
In that case, a fallback is necessary.
A safe strategy here is to start continuation with a relatively small prediction distance and increase it by a factor typically smaller than $1.4$ until the number of Newton iterations necessary for correction increases to a target number past which the algorithm is thought inefficient.
This number is usually $3$ or $4$.
If the number of Newton iterations becomes larger, the prediction distance is multiplied by a factor around $0.7$ to make sure that the simulation always runs simple incremental but non trivial calculations.

\subsubsection{Pseudo-arclength continuation}

Pseudo arc-length continuation does not work with fixed parameter values but rather adds a scalar equation to system (\ref{fu0}), making it an $(n+1) \times (n+1)$ problem.
The basic idea is to allow $\lambda$ to vary while imposing a condition on the arclength parameter $\triangle s$ \cite{Keller77,Seydel09}.
Again, several types of such conditions can be chosen, among which the orthogonal correction:
\begin{equation}
\label{withtang}
\dot{\bf u}_i \cdot ({\bf u}^j - {\bf u}_i) + \dot{\lambda} (\lambda^j - \lambda_i) = \triangle s,
\end{equation}
where $({\bf u}_i,\lambda_i)$ is the last known solution along the branch and $(\dot{\bf u}_i,\dot{\lambda})$ a normalized expression of the tangent to the branch at this point.
This condition imposes that the correction be orthogonal to the tangent to the branch at a distance $\triangle s$ from $({\bf u}_i,\lambda_i)$.

One way to approximate $(\dot{\bf u}_i,\dot{\lambda})$ is to solve the linear system:
\begin{equation}
\label{tangent}
\begin{pmatrix}
D_{{\bf u}} {\bf F} ({\bf u}_i,\lambda_i) & D_{\lambda} {\bf F} ({\bf u}_i,\lambda_i)
\end{pmatrix}
\begin{pmatrix}
\dot{\bf u}_i \\
\dot \lambda_i
\end{pmatrix}
=
{\bf 0}
,
\end{equation}
where the derivatives are evaluated at the point $({\bf u}_i,\lambda_i)$ and where $|(\dot{\bf u}_i,\dot{\lambda})| = 1$.
This implies an additional step of linear algebra and might yield significant additional computing time for large scale systems.
Instead, the tangent can be approximated by interpolation from previously calculated points along the branch \cite{Dijkstra14}.
This accuracy sacrifice does not hinder continuation.
It represents a less elegant but more efficient choice.

Once the tangential vector is evaluated, the $(n+1) \times (n+1)$ system of the $j$-th Newton iteration reads:
\begin{equation}
\begin{pmatrix}
D_{{\bf u}} {\bf F} ({\bf u}^j,\lambda^j) & D_{\lambda} {\bf F} ({\bf u}^j,\lambda^j) \\
\dot{\bf u}_i & \dot{\lambda}
\end{pmatrix}
\begin{pmatrix}
\delta {\bf u}^j \\
\delta \lambda^j
\end{pmatrix}
=
\begin{pmatrix}
{\bf F} ({\bf u}^j,\lambda^j) \\
\sigma
\end{pmatrix}
,
\end{equation}
where $\sigma = \dot{\bf u}_i \cdot ({\bf u}^j - {\bf u}_i) + \dot{\lambda} (\lambda^j - \lambda_i) - \triangle s$ is the condition residual.
Lastly, the correction is applied: ${\bf u}^{j+1} = {\bf u}^j - \delta {\bf u}^j$ and $\lambda^{j+1} = \lambda^j - \delta \lambda^j$.
A sketch of the way pseudo-arclength continuation works is shown in figure \ref{continuation}(b) alongside fixed parameter continuation to highlight the differences.

The advantage of pseudo-arclength continuation over fixed parameter continuation is that it parametrizes directly an approximation of $\mathcal{B}(s)$ so does not fail at special points.
On the other hand, to retain the whole essence of the method, the prediction has to involve a condition on the arclenth like equation (\ref{arcpredic2}) or the application of equation (\ref{withtang}) to the prediction $({\bf u}_i^p,\lambda_i^p)$ instead of the running Newton iteration.
The drawback of this is the creation of additional linear systems to be solved and thus to a substantial increase in computing time, especially in high dimension algebra.
Note that the increase in the dimension of the pseudo-arclength algorithm compared to the fixed parameter one from $n \times n$ to $(n+1) \times (n+1)$ has a negligible impact in terms of performance on such systems.

%%%%%%%%%%%%%%%%%%%%%%%%%%%%%%%%%%%%%%%%%%%%%%%%%%%%%%%%%%%%%%%%%%%%%%%
%%%%%%%%%%%%%%%%%%%%%%%%%%%%%%%%%%%%%%%%%%%%%%%%%%%%%%%%%%%%%%%%%%%%%%%
\section{The preconditioning method}\label{precon}
%%%%%%%%%%%%%%%%%%%%%%%%%%%%%%%%%%%%%%%%%%%%%%%%%%%%%%%%%%%%%%%%%%%%%%%
%%%%%%%%%%%%%%%%%%%%%%%%%%%%%%%%%%%%%%%%%%%%%%%%%%%%%%%%%%%%%%%%%%%%%%%

In general, incompressible fluid flows are modeled using the Navier--Stokes equation together with the continuity equation:
\begin{eqnarray}
\label{ns1}
&\partial_t {\bf u} + \left( {\bf u} \cdot \nabla \right) {\bf u} = - \nabla p + \frac{1}{Re} \nabla^2 {\bf u},\\
\label{ns2}
&\nabla \cdot {\bf u} = 0,
\end{eqnarray}
where $t$ represents time, ${\bf u}$ the velocity field, $p$ the pressure and $Re$ the dimensionless Reynolds number.
The choice of nondimensionalization is here purely informative and does not impact the method.
The discretized version of these equations yields large dynamical systems, often exceeding $10^5$ degrees of freedom which are typically solved using the Newton iteration.
It follows that iterative methods are preferred for the inversion of the Jacobian.
In fact, evaluating the Jacobian itself would require memory space that is hardly affordable\footnote{The storage of a $10^5 \times 10^5$ double precision matrix takes about $80$GB.} and alternative methods that do not necessitate the evaluation of the Jacobian can be used.
These methods are called {\it matrix free methods} and revolve around the ability to express the product between the Jacobian and an arbitrary vector without ever evaluating the Jacobian itself.

Another issue arises when dealing with incompressible flows: the Jacobian is often ill-conditioned.
This feature is in fact common to flows that are highly diffusive, i.e., for which the Reynolds number $Re$ is small and the Laplacian term is dominant.
An efficient preconditioner for such systems is the Stokes preconditioner \cite{Tuckerman89,Mamun95}.
This preconditioner is naturally implementable within a matrix free method which makes it a method of choice for continuation of diffusion-dominated incompressible flows.
For problems that are not dominated by diffusion, however, another preconditioner has to be used.

We consider the following dynamical system:
\begin{equation}
\label{stoeq}
\partial_t {\bf u} = N({\bf u}) + L {\bf u},
\end{equation}
where $t$ is time, ${\bf u}$ is the solution field, $N({\bf u})$ represents a nonlinear term and $L {\bf u}$ a linear term.

We use the first order implicit Euler scheme:
\begin{equation}
\frac{{\bf u}^{t+\triangle t} - {\bf u}^{t}}{\triangle t} = N({\bf u}^{t}) + L {\bf u}^{t+\triangle t},
\end{equation}
where $\triangle t$ is the timestep and ${\bf u}^{t}$ is the evaluation of ${\bf u}$ at time $t$.
On expressing ${\bf u}^{t+\triangle t}$, we get:
\begin{equation}
{\bf u}^{t+\triangle t} = \left( I - \triangle t \, L \right)^{-1} \left[ {\bf u}^{t} + \triangle t \, N({\bf u}^{t}) \right].
\end{equation}
The preconditioner is obtained by substracting ${\bf u}^{t}$ from ${\bf u}^{t+\triangle t}$:
\begin{eqnarray}
&{\bf u}^{t+\triangle t} - {\bf u}^{t} &= \left( I - \triangle t \, L \right)^{-1} \left[ {\bf u}^{t} + \triangle t \, N({\bf u}^{t}) \right] - {\bf u}^{t}\\
& &= \left( I - \triangle t \, L \right)^{-1} \left[ {\bf u}^{t} + \triangle t \, N({\bf u}^{t}) - {\bf u}^{t} + \triangle t \, L {\bf u}^{t} \right]\\
\label{sto1}
& &= \triangle t \left( I - \triangle t \, L \right)^{-1} \left[ N({\bf u}^{t}) + L {\bf u}^{t} \right],
\end{eqnarray}
where the right-hand-side is the evaluation at time $t$ of the right-hand-side of equation (\ref{stoeq}).

The general form for the preconditioned equation is thus:
\begin{equation}
\label{adaptprecond}
{\bf u}^{t + \triangle t} - {\bf u}^{t} = c \, P^{-1} \left[ N({\bf u}^{t}) + L {\bf u}^{t} \right],
\end{equation}
where $c$ is a constant and $P$ the preconditioner.
The values taken by these quantities are summarized in table \ref{tabprecon} as $\triangle t$ is changed. 
%%%%%%%%%%%%%%%%%%%%%%%%%%%%%%%%%%%%%%%%%%%%%%%%%%%
\begin{table}
\caption{Preconditioner $P$ and constant $c$ for the two limits and intermediate values of $\triangle t$.}
\begin{center}\footnotesize
\renewcommand{\arraystretch}{1.3}
\begin{tabular}{|P{3cm}|P{6cm}|P{3cm}|} \hline
\multicolumn{3}{|c|}{${\bf u}^{t + \triangle t} - {\bf u}^{t} = c P^{-1} \left[ N({\bf u}^{t}) + L {\bf u}^{t} \right]$}\\ \hline
$\triangle t \ll 1$ & $\triangle t = O(1)$ & $\triangle t \gg 1$ \\ \hline
$c=\triangle t$ & $c=\triangle t$ & $c=1$ \\ \hline
$P \rightarrow I$ & $P=I-\triangle t L$ & $P \rightarrow -L$ \\ \hline
\end{tabular}
\end{center}
\label{tabprecon}
\end{table}
%%%%%%%%%%%%%%%%%%%%%%%%%%%%%%%%%%%%%%%%%%%%%%%%%%%

There are two asymptotic regimes.
First, for sufficiently small $\triangle t$, $P$ approaches the identity operator and the equation is solved without effective preconditioning.
Second, assuming $\triangle t \gg 1$, it follows $\left( I - \triangle t L \right)^{-1} \approx - \left( \triangle t L \right)^{-1}$ and expression (\ref{sto1}) becomes:
\begin{equation}
\label{rhssto}
{\bf u}^{t+\triangle t} - {\bf u}^{t} \approx - L^{-1} \left[ N({\bf u}^{t}) + L {\bf u}^{t} \right],
\end{equation}
thus providing a Laplacian preconditioner.
This limit is known as Stokes preconditioner \cite{Tuckerman89,Mamun95} and was first used by Mamun \& Tuckerman to study symmetry breaking instabilities in spherical Couette flow and proved efficient even at relatively large Reynolds numbers\footnote{The authors computed solutions up to $Re = 2200$.} \cite{Mamun95}.
Despite the development of such a method for shear flows, rare are the subsequent uses of the Stokes preconditioner in this field.
Some studies have however emerged, like the investigation of the three-dimensional instability of a flow passed a step \cite{Barkley02}.
The Stokes preconditioner has on the other hand been used extensively to study pattern formation in convection where the Reynolds number is often unity, due to other choices of dimensionalization \cite{Bergeon98,Bergeon02,Batiste06,Mercader09,Boronska10,Beaume11,Feudel11,Lojacono11,Mercader11,Beaume13ddd1,Beaume13ddd2,Beaume13rc1,Beaume13rc2,Torres13,Torres14}.

Previous work on weakly diffusive shear flows have not necessitated preconditioning \cite{Wang07}.
In intermediate regimes, however, difficulties may arise: diffusion may not be weak enough for the system to be solved successfully without preconditioner, but the Stokes preconditioner overcompensates and yields an ill-conditioned system again.
To address this issue, we extend the Stokes preconditioner to non-asymptotic values of $\triangle t$.
The resulting preconditioner is $P = \left( I - \triangle t L \right)$, with $c=\triangle t$ and proves efficient for a wide range of systems, from coupled convection to shear flows, as we shall see.

Similarly, the linearization of equation (\ref{stoeq}),
\begin{equation}
\label{stoeqlin}
\partial_t \delta{\bf u} = \delta N({\bf u}) \, \delta {\bf u}+ L \delta{\bf u},
\end{equation}
where $\delta N({\bf u})$ is the linearization of the nonlinear term evaluated at ${\bf u}$ and $\delta {\bf u}$ an elementary displacement in ${\bf u}$, can be treated the same way to yield
\begin{eqnarray}
&\delta{\bf u}^{t+\triangle t} - \delta{\bf u}^{t} &\approx c \, P^{-1} \left[ \delta N({\bf u}^{t}) \, \delta {\bf u}^{t}+ L \delta{\bf u}^{t} \right]\\
\label{jacobsto}
& &\approx c \, P^{-1} J({\bf u}^{t}) \, \delta{\bf u}^{t},
\end{eqnarray}
where $J({\bf u}^{t}) = \delta N({\bf u}^{t}) + L$ is the Jacobian of system (\ref{stoeq}).
We can then compute the preconditioned Jacobian of equation (\ref{stoeq}) by computing one step forward in time of the linearized equation (\ref{stoeqlin}) using an implicit Euler scheme and then substracting the initial condition.
The preconditioner used for the base and linearized equations is the same for the same value of $\triangle t$.

To search for stationary flow solutions, we consider the following Newton method:
\begin{equation}
J({\bf u}) \, \delta {\bf u} \approx N({\bf u}) + L {\bf u},
\end{equation}
with correction ${\bf u} = {\bf u} - \delta {\bf u}$.
By multiplying this equation by $c \, P^{-1}$, we obtain:
\begin{equation}
\label{stofinal}
c \, P^{-1} J({\bf u}) \, \delta {\bf u} \approx c \, P^{-1} \left( N({\bf u}) + L {\bf u} \right).
\end{equation}
The left-hand-side of equation (\ref{stofinal}) is evaluated using the implicit Euler scheme given by equation (\ref{jacobsto}) and the right-hand-side of equation (\ref{stofinal}) is computed using equation (\ref{adaptprecond}).
The inversion of the Jacobian can be performed using iterative methods \cite{Saad03} such as the biconjugate gradient stabilized method \cite{Vandervorst92}.

%%%%%%%%%%%%%%%%%%%%%%%%%%%%%%%%%%%%%%%%%%%%%%%%%%%%%%%%%%%%%%%%%%%%%%%
%%%%%%%%%%%%%%%%%%%%%%%%%%%%%%%%%%%%%%%%%%%%%%%%%%%%%%%%%%%%%%%%%%%%%%%
\section{Examples}\label{secex}
%%%%%%%%%%%%%%%%%%%%%%%%%%%%%%%%%%%%%%%%%%%%%%%%%%%%%%%%%%%%%%%%%%%%%%%
%%%%%%%%%%%%%%%%%%%%%%%%%%%%%%%%%%%%%%%%%%%%%%%%%%%%%%%%%%%%%%%%%%%%%%%

In this section, two examples of continuation of incompressible fluid flows are considered.
Numerical continuation is performed using the preconditioner presented in section \ref{precon}.
The first flow considered is doubly diffusive convection, in which the preconditioner is tested in a three-dimensional configuration.
The second test flow is a two-dimensional model of shear flow in which a parameter dependent use of the preconditioner is prescribed and for which the system of equations is split into two sets preconditioned in different manners but solved simultaneously.

%%%%%%%%%%%%%%%%%%%%%%%%%%%%%%%%%%%%%%%%%%%%%%%%%%%%%%%%%%%%%%%%%%%%%%%
%%%%%%%%%%%%%%%%%%%%%%%%%%%%%%%%%%%%%%%%%%%%%%%%%%%%%%%%%%%%%%%%%%%%%%%
\subsection{Doubly diffusive convection}
%%%%%%%%%%%%%%%%%%%%%%%%%%%%%%%%%%%%%%%%%%%%%%%%%%%%%%%%%%%%%%%%%%%%%%%
%%%%%%%%%%%%%%%%%%%%%%%%%%%%%%%%%%%%%%%%%%%%%%%%%%%%%%%%%%%%%%%%%%%%%%%

We consider a Boussinesq fluid constituted of two components, the heaviest of which is referred to as a salt.
The fluid is placed within a three-dimensional enclosure of square horizontal cross-section and aspect ratio $19.8536$ in the vertical direction.
The flow is driven by buoyancy through the imposition of large scale horizontal gradients of temperature and concentration: one wall is maintained at a larger temperature and salinity than the opposite one.
The other walls are modelled using no flux conditions.
No-slip boundary conditions are imposed at all walls.

The nondimensional equations governing the dynamics of this flow are:
\begin{eqnarray}
&Pr^{-1} \left[ \partial_t {\bf u} + ({\bf u} \cdot \nabla){\bf u} \right] = -\nabla p + Ra (T-C) {\bf \hat{x}} + \nabla^2 {\bf u},\\
&\nabla \cdot {\bf u} =0,\\
&\partial_t T + ({\bf u} \cdot \nabla)T = \nabla^2 T,\\
&\partial_t C + ({\bf u} \cdot \nabla)C = \tau \nabla^2 C,
\end{eqnarray}
where $t$ is time, ${\bf u} \sim (u,v,w)$ is the velocity field in the Cartesian frame $({\bf \hat{x}},{\bf \hat{y}},{\bf \hat{z}})$, $p$ is the pressure, and $T$ and $C$ are linear rescaling of the fluid temperature and salt concentration in the Boussinesq approximation.
Here, ${\bf \hat{x}}$ represents the vertical unit vector in the ascending direction.
In addition to these quantities, three nondimensional parameters are introduced.
The Prandtl number $Pr$ is the ratio of the kinematic viscosity over the thermal diffusivity, the inverse Lewis number $\tau$ is the ratio of the salt diffusivity over the thermal diffusivity and the Rayleigh number $Ra$ quantifies the buoyancy strength and will be used as the continuation parameter in what follows.
These equations are complemented with boundary conditions:
\begin{eqnarray}
\label{ddbc1}
&{\rm at}~x=\{0,L\}~{\rm or}~y=\{0,1\}&: u = v = w = \partial_n T=\partial_n C = 0,\\
\label{ddbc3}
&{\rm at}~z=0&: u = v = w = T = C = 0,\\
&{\rm at}~z=1&: u = v = w = T -1 = C -1 = 0,
\end{eqnarray}
where the operator $\partial_n$ represents the spatial derivative in the direction normal to the wall.
More details on the physical setup are available in \cite{Beaume13ddd1}.

\subsubsection{Numerics}

\begin{figure}
\begin{center}
\includegraphics[width=0.8\textwidth]{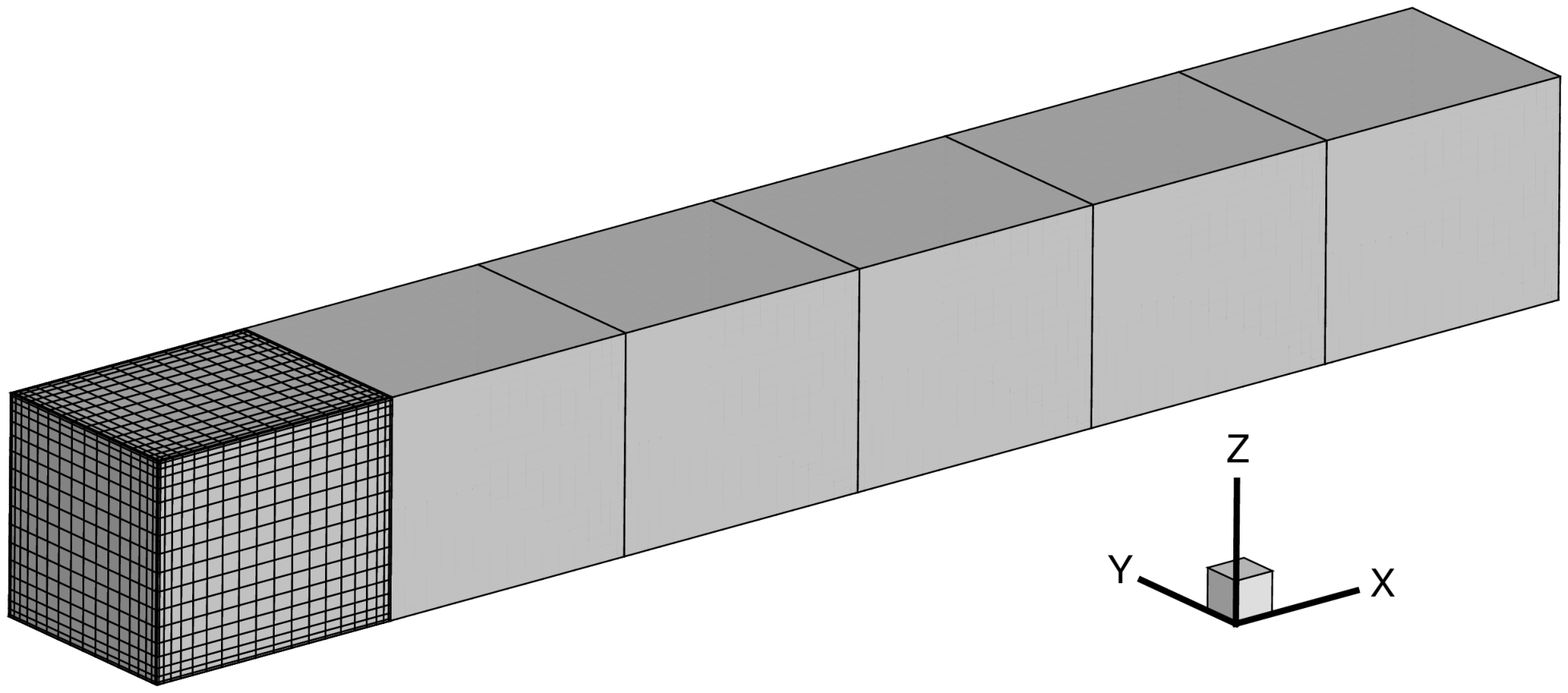}
\end{center}
\caption{Example of discretization for the doubly diffusive convection problem. For the sake of the representation, only $6$ out of the $16$ elements are shown with only the first (bottom left) element displaying its surface mesh. The inner grid is hidden to avoid overloading the figure.}
\label{mesh}
\end{figure}

The physical domain is meshed using $16$ identical spectral elements of size $l_x\approx 1.24$, $l_y=l_z=1$.
Each element is meshed using Gauss--Lobato--Legendre points in all three direction: $21$ in $x$, $19$ in $y$ and $z$.
The discretization strategy is illustrated in figure \ref{mesh} and yields $16 \times 21 \times 19 \times 19 \times 6 = 727,776$ degrees of freedom (counting in the pressure).
%%%%%%%%%%%%%%%%%%%%%%%%%%
%%%%%%%%%%%%%%%%%%%%%%%%%%

Time is discretized using a straightforward first order Euler scheme for temperature and concentration:
\begin{eqnarray}
&T^{(n)} = \left( I - \triangle t \, \nabla^2 \right)^{-1} \left(T^{(n-1)} - \triangle t \, [ ({\bf u} \cdot \nabla )T]^{(n-1)} \right),\\
&C^{(n)} = \left( I - \triangle t \, \tau \, \nabla^2 \right)^{-1} \left(C^{(n-1)} - \triangle t \, [ ({\bf u} \cdot \nabla )C]^{(n-1)} \right),
\end{eqnarray}
where $T^{(n)}$ stands for the evaluation of the temperature at the $n$-th timestep, $I$ is the identity operator and $\triangle t$ is the timestep.

The incompressible Navier--Stokes equation is discretized using a first order splitting method described by Karniadakis, Israeli and Orszag \cite{Karniadakis91}.
An intermediate velocity is predicted which takes into account buoyancy and advection:
\begin{equation}
{\bf \hat{u}} = {\bf u}^{(n-1)} - \triangle t \, [ ({\bf u} \cdot \nabla ){\bf u}]^{(n-1)} + \triangle t \, Ra \, (T-C)^{(n-1)} {\bf \hat{x}},
\end{equation}
which is corrected using the incompressibility condition by introducing the velocity ${\bf \hat{\hat{u}}}$:
\begin{equation}
\label{kio2}
{\bf \hat{\hat{u}}} = {\bf \hat{u}} - \triangle t \, \nabla p^{(n)},
\end{equation}
where the pressure is defined by the Poisson problem obtained by taking the divergence of equation (\ref{kio2}):
\begin{equation}
\nabla^2 p^{(n)} = \frac{1}{\triangle t} \nabla \cdot {\bf \hat{u}},
\end{equation}
complemented with the boundary condition:
\begin{equation}
\partial_n p^{(n)} = \left( (T-C)^{(n-1)} {\bf \hat{x}} - [ ({\bf u} \cdot \nabla ){\bf u}]^{(n-1)} - \nabla \times \nabla \times {\bf u}^{(n-1)} \right) \cdot {\bf \hat{n}},
\end{equation}
where, ${\bf \hat{n}}$ represents the vector normal to the boundary and the last term is the reduction of the Laplacian term using the divergence condition.
The time-step is completed via the following operation:
\begin{equation}
{\bf u}^{(n)} = \left( I - \triangle t \, \nabla^2 \right)^{-1} (\triangle t {\bf \hat{\hat{u}}} ),
\end{equation}
where the original boundary conditions (\ref{ddbc1})--(\ref{ddbc3}) are used.

Time-stepping then only requires the inversion of Helmholtz operators coming from the spatial discretization.
As a consequence of the choice of spectral elements, the Helmholtz operators are sparse tensors and a Schur decomposition has been performed to invert them efficiently.
Continuation is performed on the solution vector $({\bf u},T,C) = (u,v,w,T,C)$ using the temporal schemes above and the same $\triangle t$ for all equations.

\subsubsection{Results}

This flow configuration exhibits localized pattern formation at onset through a subcritical bifurcation.
The branches emerging from this bifurcation produce well-bounded back and forth oscillations in parameter space in a behavior known as {\it snaking}.
More detailed information is available in \cite{Beaume13ddd1}.

\begin{figure}
\begin{center}
\includegraphics[width=\textwidth]{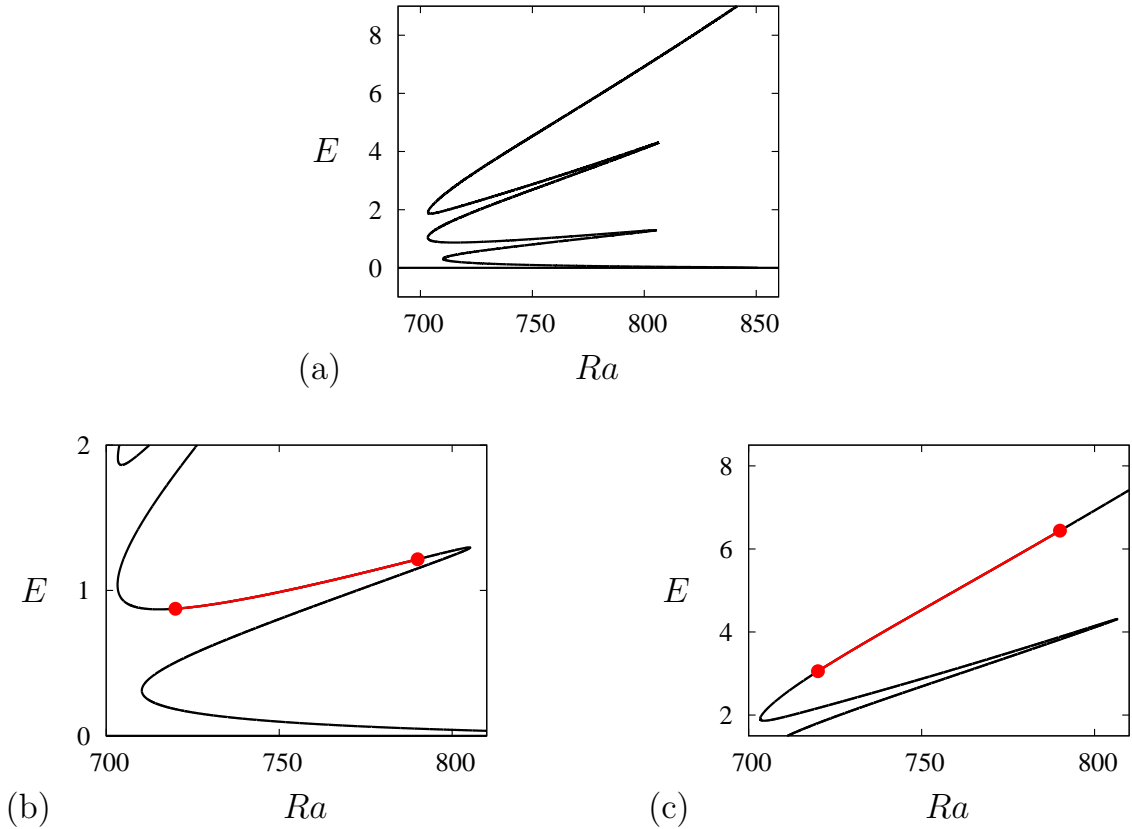}
\end{center}
\caption{(a) Bifurcation diagram showing the kinetic enery $E$ of the trivial motionless fluid branch (horizontal line) and of one of the localized state branches as a function of the Rayleigh number $Ra$. (b) Zoom in the first region of interest: the test of the preconditioner is run through leftward continuation along the red portion of the branch delimited by the dots. (c) Zoom in the second region of interest where the tests is run rightwards. The solutions denoted by a dot are shown in figure \ref{locasols}.}
\label{locabranch}
\end{figure}

\begin{figure}
\begin{center}
\includegraphics[width=\textwidth]{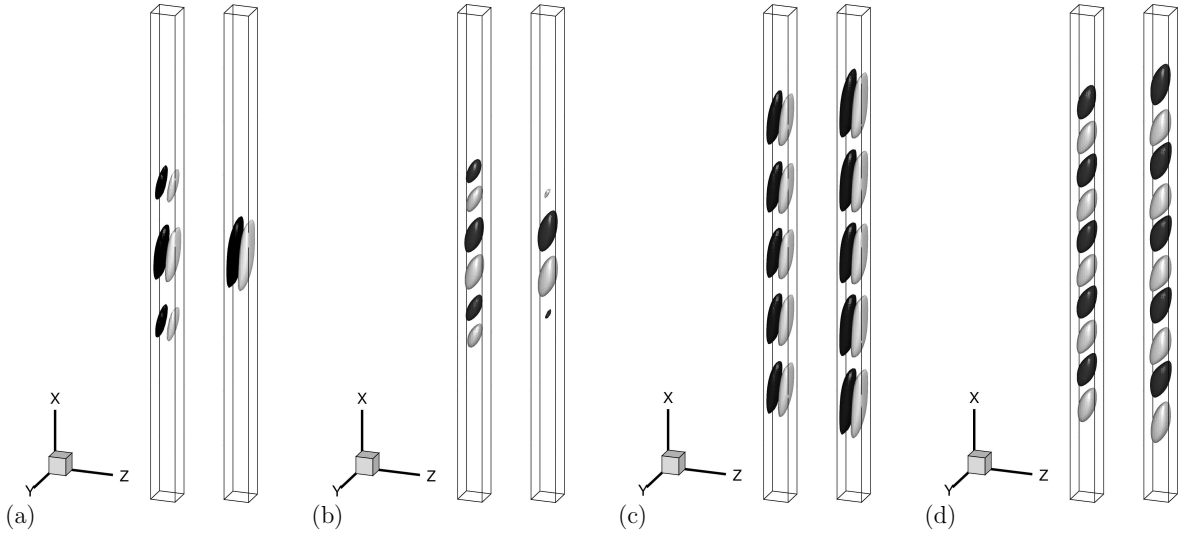}
\end{center}
\caption{Isosurfaces of the $x$-velocity $u=\pm0.5$ (a,c) and of the $z$-velocity $w=\pm0.2$ (b,d). In both cases, the solutions represent the extrema of the test segments shown in figure \ref{locabranch}(b) and (c) with the left solution being taken at $Ra \approx 720$ and the right solution at $Ra \approx 790$. Panels (a) and (b) represent the end points of the test along the lower segment (figure \ref{locabranch}(b)), while panels (c) and (d) represent those along the upper segment (figure \ref{locabranch}(c)). The light (resp. dark) color indicates the positive (resp. negative) contour.}
\label{locasols}
\end{figure}

We focus here on one branch of spatially localized states.
The branch is shown in figure \ref{locabranch}.
%%%%%%%%%%%%%%%%%%%%%%%%%%
%%%%%%%%%%%%%%%%%%%%%%%%%%
To design a test for the preconditioner, we select two segments along the branch and impose fixed parameter continuation to the algorithm to avoid intricacies related to continuation with a variable parameter.
One of this segment consists in solutions that are very localized and where most of the domain is filled with ``zeros'' while the other segment consists in domain-fillin solutions.
The saddle-nodes occur at $Ra \approx 703$ and $Ra \approx 807$ so we restrict continuation to $720 < Ra < 790$.
The solution at $Ra \approx 790$ and at $Ra \approx 720$ for both segments are shown in figure \ref{locasols}.
%%%%%%%%%%%%%%%%%%%%%%%%%%
%%%%%%%%%%%%%%%%%%%%%%%%%%
These solutions display a number of similar convection rolls centered in the domain.
In each of these rolls, the fluid goes up along the hot and saltier wall at $z=1$ and down on the opposite wall.
The flow is not purely two-dimensional: due to the presence of walls at $y=0$ and $y=1$, a weak flow in the $y$ direction is generated.
This flow is typically an order of magnitude lower than in the two other direction and is therefore not shown here.
As the branch is continued along the lower segment, the Rayleigh number is decreased and the left the solution changes from the right panels to the left panels of figure \ref{locasols}, thereby adding rolls on either side of the central roll.
Along the upper segment, the solution is continued in the direction of increasing Rayleigh numbers (from the left to the right panel in figure \ref{locasols}) and the rolls grow in size and amplitude.

The algorithmic parameters are kept at the values used during the original study \cite{Beaume13ddd1} and which were determined using a combination of intuition and parametric benchmark: the tolerance of the BCGStab is fixed at $10^{-2}$, the continuation step is initialized at $\delta Ra = 10^{-3}$ and can go up to $1$.
The convergence of the solution is assessed by calculating the following quantity:
\begin{equation}
\mathcal{L} = \left|\frac{1 + \triangle t}{I - \triangle t \, L} \; \left[N({\bf u}) + L {\bf u} \right]\right|_{L_2},
\end{equation}
computed by stepping forward once in time (see equation (\ref{adaptprecond})) and multiplying by $(1 + \triangle t)/\triangle t$ and then taking the $L_2$-norm of the resulting vector.
The multiplying step is taken so that the relative convergence of $N({\bf u}) + L {\bf u}$ does not depend strongly on the value of $\triangle t$.
The continuation is accelerated by a factor of $1.2$ after each successful step.
If more than $4$ Newton iterations were necessary to obtain $L < 10^{-7}$, the continuation step remains untouched.
If the Newton iteration fails to converge or if the number of gradient iteration needed to invert the Jacobian is greater than $5000$, the current step is cancelled and another attempt with a continuation step $10\%$ smaller is made. 

\begin{figure}
\begin{center}
\includegraphics[width=\textwidth]{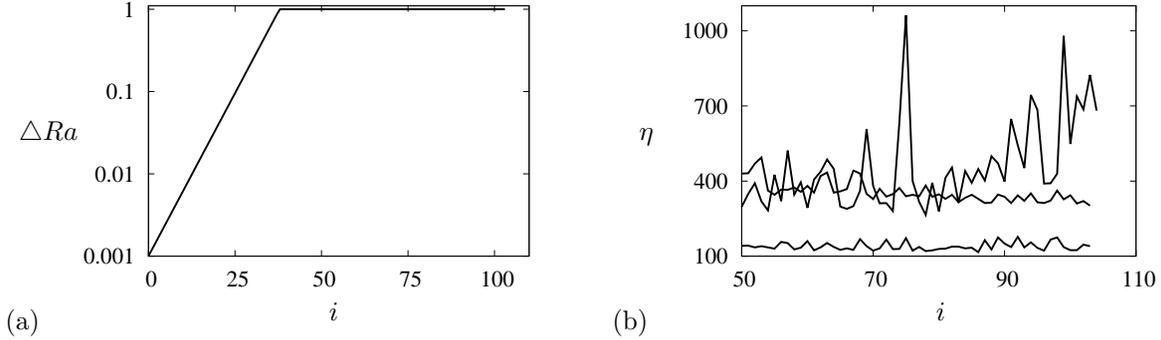}
\end{center}
\caption{(a) Evolution of the parameter step $\triangle Ra$ as a function of the continuation point $i$ for $\triangle t = 0.06$. It is exponential with slope $1.2$ until $\triangle t = 1$ is reached where the growth is algorithmically frozen. Most of the other simulations yielded the same results. (b) Number of gradient iterations $\eta$ needed to converge as a function of the continuation point $i$ for $i \ge 50$ and for $\triangle t = 0.003$ (upper curve displaying large variations), $\triangle t = 0.06$ (lower curve) and $\triangle t = 10^5$ (upper curve with small variations). These results have been obtained during continuation of the lower segment from figure \ref{locabranch}(b).}
\label{dddinitmax1}
\end{figure}

A number of attempts were made for both continuation segments, with $\triangle t$ ranging from $10^{-4}$ to $10^{8}$.
The smallest values of $\triangle t$ did not allow the algorithm to converge, implying that this problem does indeed need preconditioning.
I report here the successful simulations, for $\triangle t \ge 10^{-3}$.
The basic algorithmic behavior is illustrated in figure \ref{dddinitmax1} on simulations carried out on the lower segment in figure \ref{locabranch}(b).
%%%%%%%%%%%%%%%%%%%%%%%%%%
%%%%%%%%%%%%%%%%%%%%%%%%%%
Panel (a) shows the continuation step as a function of the iteration number $i$ and indicates an exponential acceleration until $\triangle Ra = 1$.
Nearly all the simulations run provided the same results here, which indicates that the continuation neither fails nor become marginally successful in these conditions.
The few exceptions were generally obtained for $\triangle t \le 0.01$ for which the preconditioner is not strong.
These cases display a few continuation points at which more than $4$ Newton iterations are needed and where therefore $\triangle Ra$ remains fixed before continuing its progression to $1$.
Figure \ref{dddinitmax1}(b) shows for three cases the typical values taken by $\eta$, the total number of gradient iterations needed to converge a continuation point (summed up on all the Newton steps required at this continuation point).
For $\triangle t = 0.003$, $39$ continuation points are necessary to reach $\triangle Ra = 1$.
The convergence speed displays large variations, with $\eta$ varying between $265$ and $1062$  and averaging $\bar{\eta} \approx 471$ with a standard deviation of $\sigma \approx 179$ over the last $50$ continuation points.
Results at $\triangle t = 10^5$ do not suffer from such large oscillations, $\sigma \approx 32$, and their average number of iteration is significantly lower: $\bar{\eta} \approx 349$.
The best results were obtained for $\triangle t = 0.06$ (lower curve in figure \ref{dddinitmax1}(b)): $\bar{\eta} \approx 140$ with $\sigma \approx 17$.
 
\begin{figure}
\begin{center}
\includegraphics[width=\textwidth]{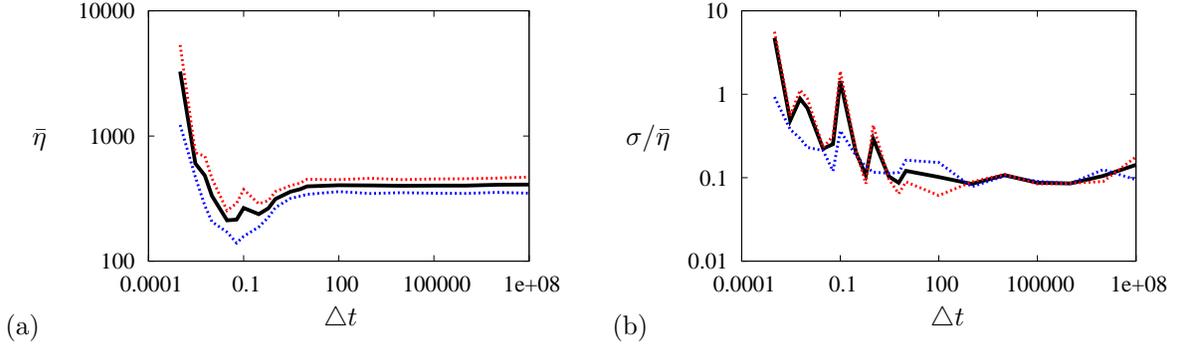}
\end{center}
\caption{Compilation of the results showing the average number of gradient iterations needed to converge the Newton method on a sample consisting of the last $50$ continuation points $\bar{\eta}$ as a function of $\triangle t$ (a). The right panel (b) shows the standard deviation $\sigma$ associated with these simulations normalized by $\bar{\eta}$ as a function of $\triangle t$. In both panels, the red dashed curve corresponds to the upper continuation segment results while the blue dashed curve corresponds to the lower continuation segment results. The thick black curve indicates the average between both segments.}
\label{dddresmax1}
\end{figure}

The results of all the simulations are compiled in figure \ref{dddresmax1} and presented through the average number of gradient iterations needed to converge one continuation point $\bar{\eta}$ and the normalized standard deviation from this result $\sigma/\bar{\eta}$ for the last $50$ continuation points, i.e., when the continuation step is constant: $\triangle Ra = 1$.
%%%%%%%%%%%%%%%%%%%%%%%%%%
%%%%%%%%%%%%%%%%%%%%%%%%%%
Figure \ref{dddresmax1}(a) clearly indicates the presence of three distinct regions where the preconditioner behaves differently.
The first and most obvious is obtained for $\triangle t < 0.01$: for such low values of $\triangle t$, there is effectively little preconditioning in place and on such a diffusive system, the continuation struggles as shown by the large increase of $\hat{\eta}$ as $\triangle t$ decreases.
For $\triangle t \le 0.001$, continuation becomes impossible, confirming the need for a preconditioner.
A second region of interest is obtained at large $\triangle t$ for which the preconditioner takes the form of a Laplacian.
In fact, little difference is observed for simulations run with $\triangle t \ge 10^2$, indicating that the asymptotic regime of the Laplacian preconditioner is reached.
In this regime, converging a continuation point costs around $350$ gradient iterations for the lower segment and $455$ for the upper one.
Between these two regions lies a sweet spot located around $\triangle t = 0.06$.
This sweet spot consists in an interval where the preconditioner is at its performance peak.
The most efficient continuation of the lower segment was performed at $\triangle t = 0.06$ with an average of $\hat{\eta} \approx 140$ gradient iterations per continuation point.
Similarly, the most efficient continuation of the upper segment corresponded to $\triangle t = 0.03$ and $\hat{\eta} \approx 253$.

Furthermore, figure \ref{dddresmax1}(b) shows the standard deviation associated to the number of gradient iterations normalized by the average to quantify the variation significance.
The Laplacian preconditionner is very robust with a relative standard deviation of about $10\%$ of $\hat{\eta}$ for $\triangle t \ge 10^2$.
As $\triangle t$ decreases, $\sigma / \hat{\eta}$ increases to reach values above $1$ for the smallest successful $\triangle t$.
This highlights one particular characteristic: although the preconditioner is at its best in shear performance for values of $\triangle t$ around $0.06$, it is less robust there than it is for large $\triangle t$.
This accounts for the fact that the results presented in figure \ref{dddresmax1}(a) are smooth at large $\triangle t$ but display some anomalies at lower $\triangle t$.
Notice that the anomaly reported at $\triangle t = 0.1$ along the upper segment is due to the algorithm failing to converge at one point. 
It is reported here for complete transparency on the results.

\begin{figure}
\begin{center}
\includegraphics[width=0.5\textwidth]{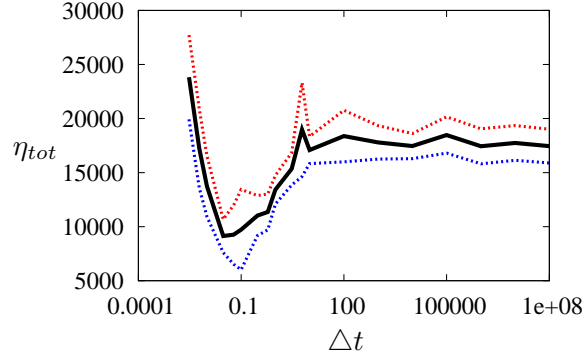}
\end{center}
\caption{Total number of gradient iterations $\eta_{tot}$ needed to complete the lower (upper) continuation segment in dashed blue (red) lines versus the algorithmic parameter $\triangle t$. The solid black line represents the average between the two data sets.}
\label{dddinfinite}
\end{figure}

Best results on a given continuation segment are obtained without bounding the continuation step.
Although this is true in principle, in practice one needs to bound it to ensure efficient continuation of saddle-nodes and to avoid jumping onto another branch in case of imperfect bifurcation.
Figure \ref{dddinfinite} shows the total number of gradient iterations undergone from the startpoint until the endpoint of the continuation segments with $\triangle t$ unbounded.
%%%%%%%%%%%%%%%%%%%%%%%%%%
%%%%%%%%%%%%%%%%%%%%%%%%%%
These results are very similar to those shown in figure \ref{dddresmax1}(a) and confirm the previous observations.

%%%%%%%%%%%%%%%%%%%%%%%%%%%%%%%%%%%%%%%%%%%%%%%%%%%%%%%%%%%%%%%%%%%%%%%
%%%%%%%%%%%%%%%%%%%%%%%%%%%%%%%%%%%%%%%%%%%%%%%%%%%%%%%%%%%%%%%%%%%%%%%
\subsection{Shear flow}
%%%%%%%%%%%%%%%%%%%%%%%%%%%%%%%%%%%%%%%%%%%%%%%%%%%%%%%%%%%%%%%%%%%%%%%
%%%%%%%%%%%%%%%%%%%%%%%%%%%%%%%%%%%%%%%%%%%%%%%%%%%%%%%%%%%%%%%%%%%%%%%

We now consider a three-dimensional fluid confined between two parallel plates of infinite extent.
The flow is driven by the imposition of a sinusoidal volume force creating shear across the fluid layer in a configuration known as {\it plane Waleffe flow}.
The Navier--Stokes equation together with the incompressibility constraint for this configuration read:
\begin{eqnarray}
\label{shearns1}
&\partial_t {\bf u} + ({\bf u} \cdot \nabla){\bf u} = -\nabla p + \frac{1}{Re} \nabla^2 {\bf u} + \frac{\sqrt{2} \pi^2}{4 Re} \sin \left(\frac{\pi y}{2}\right) {\bf \hat{x}},\\
\label{shearns2}
&\nabla \cdot {\bf u} =0,
\end{eqnarray}
where $t$ is time, ${\bf u} = (u,v,w)$ the velocity field in the $(x,y,z)$ coordinate frame where $x$ is the streamwise direction, $y$ the wall-normal direction and $z$ the spanwise direction, $p$ is the pressure and $Re$ is the Reynolds number which quantifies the imposed shear across the fluid.
These equations are accompanied with periodic boundary conditions in $x$ and $z$ and no-slip boundary conditions in $y$:
\begin{equation}
\partial_y u = v = \partial_y w = 0~{\rm at}~y=\pm 1.
\end{equation}

This flow configuration is a close cousin of plane Couette flow and is studied to investigate transition to turbulence.
A number of studies have revealed the influence of exact coherent states in the transition process \cite{Nagata90,Clever97,Waleffe03,Wang07,Schneider08,Gibson09,Kawahara12}.
These states are exact solutions of the associated system of equations and some of them follow an asymptotic behavior as the Reynolds number is increased \cite{Wang07}.
When the corresponding asymptotic expansions are applied, the three-dimensional system (\ref{shearns1}), (\ref{shearns2}) reduces down to the following two-dimensional system:
\begin{eqnarray}
\label{shearred1}
&\partial_T u_0 + J(\phi_1,u_0) = \nabla_{\perp}^2 u_0 + \displaystyle\frac{\sqrt{2}\pi^2}{4} \sin \left(\frac{\pi y}{2} \right),\\
\label{shearred2}
&\partial_T \omega_1 + J(\phi_1,\omega_1) + 2 (\partial_y^2 - \partial_z^2) (\mathcal{R}(v_1' w_1'^*)) + 2 \partial_y \partial_z (w_1' w_1'^* - v_1' v_1'^*) = \nabla_{\perp}^2 \omega_1,\\
\label{shearred3}
&\left(\alpha^2 - \nabla_{\perp}^2\right) p_1' = 2 i \alpha (v_1' \partial_y u_0 + w_1' \partial_z u_0),\\
\label{shearred4}
&\partial_t {\bf v_{1\perp}'} + i \alpha u_0 {\bf v_{1\perp}'} = - \nabla_{\perp} p_1' + \epsilon \nabla_{\perp}^2 {\bf v_{1\perp}'},
\end{eqnarray}
where $T = \epsilon t$ with $\epsilon = Re^{-1} \ll 1$, $J(\phi_1, \cdot) = \partial_y\phi_1 \partial_z \cdot - \partial_z \phi_1 \partial_y \cdot$, $\nabla_{\perp} = (\partial_y,\partial_z)$, $\nabla_{\perp}^2 = \partial_y^2 + \partial_z^2$, $\mathcal{R}(\cdot)$ indicates the real part of $\cdot$, $\cdot^*$ is the complex conjugate of $\cdot$ and $i$ is the unit imaginary number.
The fields $u_0$, $\phi_1$ and $\omega_1$ are real while ${\bf v_{1\perp}'} = (v_1',w_1')$ and $p_1'$ are complex.
To write system (\ref{shearred1})--(\ref{shearred4}), the solution has been approximated using the following asymptotics at small $\epsilon = Re^{-1}$:
\begin{eqnarray}
&u(x,y,z,t) &\sim u_0(y,z,T) + \epsilon \left(u_1(y,z,T) + u_1'(y,z,t,T) e^{i \alpha x} + c.c.\right),\\
&v(x,y,z,t) &\sim \epsilon \left(v_1(y,z,T) + v_1'(y,z,t,T) e^{i \alpha x} + c.c.\right),\\
&w(x,y,z,t) &\sim \epsilon \left(w_1(y,z,T) + w_1'(y,z,t,T) e^{i \alpha x} + c.c.\right),
\end{eqnarray}
where $u_0$, $u_1$, $v_1$ and $w_1$ are real, $u_1'$, $v_1'$ and $w_1'$ are complex, $\alpha$ is a chosen wavelength in the streamwise direction and $c.c.$ represents the complex conjugate.
The pressure has been expanded accordingly and a streamfunction $\phi_1$ and a vorticity $\omega_1$ have been introduced such that: $v_1 = -\partial_z \phi_1$, $w_1 = \partial_y \phi_1$ and $\omega_1 = \nabla_{\perp}^2 \phi_1$.
Using this asymptotic approach, the boundary conditions read:
\begin{equation}
\partial_y u_0 = \omega_1 = \phi_1 = v_1' = \partial_y w_1' = 0~{\rm at}~y=\pm 1,
\end{equation}
together with periodic boundary conditions in $z$.
For more details on the derivation, see Beaume {\it et al.} \cite{Beaume15pre}.

\subsubsection{Numerics}

The physical domain is two-dimensional and has size $L_y = 2$ and $L_z = \pi$.
It is meshed with $32$ equidistributed points in each direction and the linear operators treated using the Fast Fourier Transform in $z$ and either the Fast Cosine Transform I or the Fast Sine Transform I in $y$ depending on the boundary condition \cite{Frigo05}.
The usual 2/3 dealiasing rule is applied to prevent frequency folding.
This resulting number of unfiltered degrees of freedom is then: $32 \times 32 \times 8 = 8,192$.

The fluctuating pressure $p_1'$ is solved for as a preliminary step:
\begin{equation}
p_1'^{(n-1)} = 2 i \alpha \left(\alpha^2 - \nabla_{\perp}^2\right)^{-1} \left(v_1'^{(n-1)} \partial_y u_0^{(n-1)} + w_1'^{(n-1)} \partial_z u_0^{(n-1)}\right),
\end{equation}
where the nonlinear right-hand-side is evaluated in physical space and the linear operator inverted in frequency space.

The same first order Euler scheme as for the doubly diffusive convection problem is used to treat time dependence in the remaining equations:
\begin{eqnarray}
&u_0^{(n)} = \left( I - \epsilon \triangle t \nabla_{\perp}^2 \right)^{-1} \left[ u_0^{(n-1)} + \epsilon \triangle t \left( - J\left(\phi_1^{(n-1)},u_0^{(n-1)}\right) + \frac{\sqrt{2} \pi^2}{4} \sin \left(\frac{\pi y}{2}\right) \right) \right],\\
&\omega_1^{(n)} = \left( I - \epsilon \triangle t \nabla_{\perp}^2 \right)^{-1} \biggl[ \omega_1^{(n-1)} + \epsilon \triangle t \biggl( - J\left(\phi_1^{(n-1)},\omega_1^{(n-1)}\right) \dots\nonumber\\
&+ 2 (\partial_z^2 - \partial_y^2) (\mathcal{R}(v_1'^{(n-1)} w_1'^{(n-1)*})) + 2 \partial_y \partial_z (v_1'^{(n-1)} v_1'^{(n-1)*} - w_1'^{(n-1)} w_1'^{(n-1)*}) \biggr) \biggr],\\
&{\bf v_{1\perp}'}^{(n)} = \left( I - \epsilon \triangle t \nabla_{\perp}^2 \right)^{-1} \left( {\bf v_{1\perp}'}^{(n-1)} + \triangle t \left( - i \alpha u_0^{(n-1)} {\bf v_{1\perp}'}^{(n-1)} - \nabla_{\perp} p_1'^{(n-1)} \right) \right).
\end{eqnarray}

Continuation is carried out on the vector $(u_0,\omega_1,{\bf v_{1\perp}'}) = (u_0,\omega_1,v_1',w_1')$ but the different balance between the mean equations (\ref{shearred1}), (\ref{shearred2}) and the fluctuation equations (\ref{shearred3}), (\ref{shearred4}) implies that different preconditioners, hence different $\triangle t$ for each set of equations, were required.
Note that the use of different $\triangle t$ for different equations of the same system is mathematically correct for as long as it is solved for steady-states.
We introduce $\triangle t_1$ which is used in the preconditioning of equations (\ref{shearred1}), (\ref{shearred2}), and $\triangle t_2$ which plays a similar role for equations (\ref{shearred3}), (\ref{shearred4}) but stress out once again that continuation is done on the whole solution vector at once.

\subsubsection{Results}

\begin{figure}
\begin{center}
\includegraphics[width=0.5\textwidth]{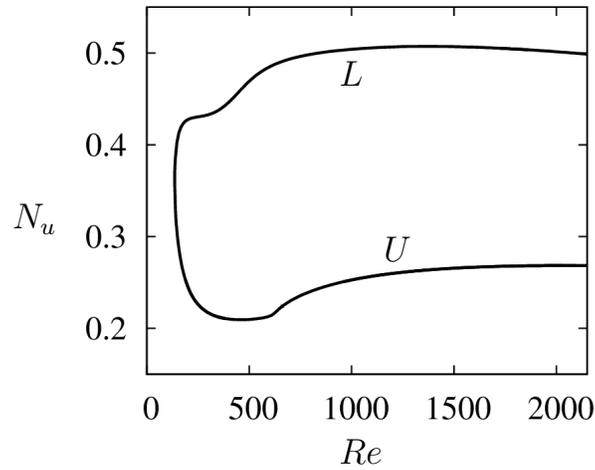}
\end{center}
\caption{Bifurcation diagram representing the double of the kinetic energy per unit volume associated with velocity $u_0$: $N_u$ versus the Reynolds number $Re$. The quantity $N_u$ is defined as follows: $N_u = D^{-1} \int_{\mathcal{D}} u_0^2 dy dz$ with $D = \int_{\mathcal{D}} dy dz$ and $\mathcal{D} = [-1;1] \times [0;L_z]$ represents the domain of integration. The trivial solution of plane Waleffe flow has $N_u =1$ (not shown). The lower (resp. upper) branch solution is labeled $L$ (resp. $U$).}
\label{shearbif}
\end{figure}

\begin{figure}
\begin{center}
\includegraphics[width=\textwidth]{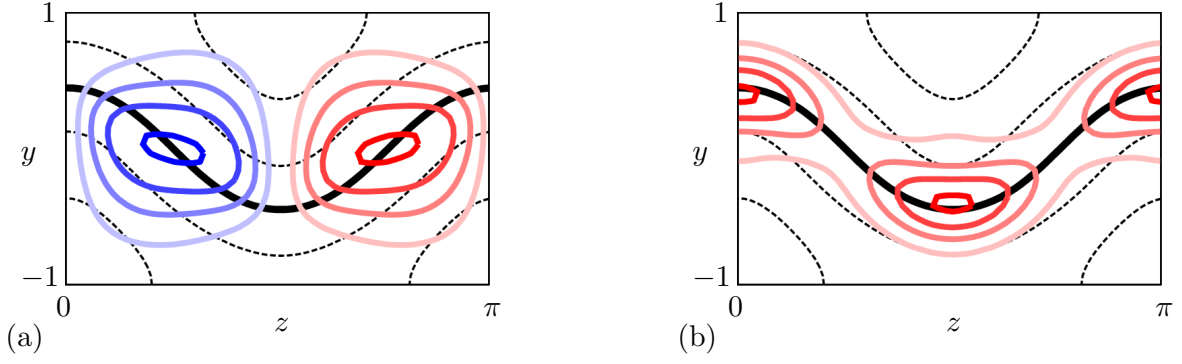}
\end{center}
\caption{Solution from the lower branch $L$ in figure \ref{shearbif} taken at $Re \approx 1000$. It is represented through equidistributed streamfunction $\phi_1$ contours (with increments of $0.4$) (a) and fluctuation amplitude $||(v_1',w_1')||$ contours (with increments of $1.75$) (b) in the $(y,z)$-plane. Positive (negative) quantities are represented in red (blue) and are plotted on top of the equidistributed contours of streamwise-invariant streamwise velocity $u_0$ (with increments of $0.5$) in black, with the thick solid black line representing the critical layer where $u_0 = 0$.}
\label{initsol}
\end{figure}

\begin{figure}
\begin{center}
\includegraphics[width=\textwidth]{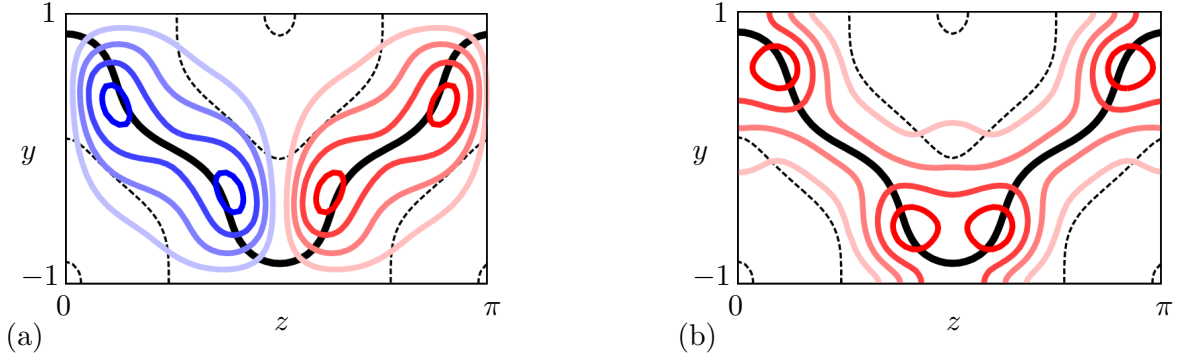}
\end{center}
\caption{Same representation as in figure \ref{initsol} but for a solution taken along the upper branch $U$ at $Re \approx 1000$. The contours are equidistributed with increments of $0.7$ for the streamfunction (a) and $2.5$ for the fluctuation amplitude (b).}
\label{shearupper}
\end{figure}

The trivial solution of plane Waleffe flow $(u_0,\omega_1,v_1',w_1') = (\frac{\sqrt{2} \pi^2}{4} \sin \left( \frac{\pi y}{2} \right),0,0,0)$ is linearly stable for all values of the Reynolds number but does not prevent other nonlinear solutions to exist.
These solutions are formed at saddle-node bifurcations at finite $Re$ and take the form of upper and lower branches, the former being energetically farther from the trivial solution than the latter.
We focus here on the most basic of these solutions obtained for a domain size of $L_x = 4 \pi$ (implying here $\alpha = 0.5$), $L_y = 2$ and $L_z = \pi$.
The bifurcation diagram is shown in figure \ref{shearbif}.
%%%%%%%%%%%%%%%%%%%%%%%%%%
%%%%%%%%%%%%%%%%%%%%%%%%%%
The solution is formed at a saddle-node at $Re \approx 136$ and splits into a lower branch state, shown in figure \ref{initsol} and an upper branch state, shown in figure \ref{shearupper}, both for $Re \approx 1000$.
%%%%%%%%%%%%%%%%%%%%%%%%%%
%%%%%%%%%%%%%%%%%%%%%%%%%%
%%%%%%%%%%%%%%%%%%%%%%%%%%
%%%%%%%%%%%%%%%%%%%%%%%%%%
As the Reynolds number is increased along these branches, the pattern remains similar but the fluctuations become sharper.

The continuation code from \cite{Beaume15pre} is modified to output relevant data for algorithmic comparison but no changes to the computational part of the code is made.
In particular, the algorithmic constants are kept the same: the tolerance for the BCGStab iterations is $10^{-2}$ with a maximum number of iteration of $1000$ before failure is declared and the convergence of the Newton method is considered reached when $\mathcal{L} < 10^{-8}$ where $\mathcal{L}$ is defined below.
%The acceleration of the continuation is set in the following way: if less than $4$ Newton iterations are necessary to converge the solution, the continuation step $\triangle Re$ is increased by $25\%$, it remains unchanged for $4$ Newton iterations and decreases by $25\%$ if more than $4$ Newton iterations were needed or if the previous continuation step has been cancelled due to the impossibility of the conjugate gradient to converge in $1000$ iterations of the Newton method to converge in $9$ iterations.
%The starting continuation step is $\triangle Re = 1$.

Due to the fundamental differences between the mean equations (\ref{shearred1}), (\ref{shearred2}) and the fluctuation equations (\ref{shearred3}), (\ref{shearred4}), different preconditioners are needed for each set.
We thus set $\triangle t_1$ and $\triangle t_2$ the respective preconditioning parameters associated to the mean and fluctuation equations.
The preconditioning method then writes:
\begin{eqnarray}
\label{shearprec1}
&u_0^{(n)} - u_0^{(n-1)} = \epsilon \, \triangle t_1 \, (I - \epsilon \, \triangle t_1 \, L_{11})^{-1} \left[L_{11} u_0^{(n-1)} + N_{11} \left(u_0^{(n-1)},\phi_1^{(n-1)}\right)\right],\\
\label{shearprec2}
&\omega_1^{(n)} - \omega_1^{(n-1)} = \epsilon \, \triangle t_1 \, (I - \epsilon \, \triangle t_1 \, L_{12})^{-1} \dots \hspace{4cm} \nonumber\\
&\hspace{4cm} \left[L_{12} \omega_1^{(n-1)} + N_{12} \left(\omega_1^{(n-1)},\phi_1^{(n-1)},{\bf v_{1\perp}'}^{(n-1)}\right)\right],\\
\label{shearprec3}
&{\bf v_{1\perp}'}^{(n)} - {\bf v_{1\perp}'}^{(n-1)} = \triangle t_2 (I - \epsilon \, \triangle t_2 \, L_2)^{-1} \dots\hspace{4cm}\nonumber\\
&\hspace{4cm} \left[\epsilon \, L_2 {\bf v_{1\perp}'}^{(n-1)} + N_2 \left({\bf v_{1\perp}'}^{(n-1)},u_0^{(n-1)},p_1'^{(n-1)}\right)\right],
\end{eqnarray}
for which 
\begin{eqnarray}
&L_{11} = L_{12} = L_2 = \nabla_{\perp}^2,\\
&N_{11} \left(u_0^{(n-1)},\phi_1^{(n-1)}\right) = - J\left(\phi_1^{(n-1)},u_0^{(n-1)}\right) + \frac{\sqrt{2} \pi^2}{4} \sin \left(\displaystyle\frac{\pi y}{2}\right),\\
&N_{12}\left(\omega_1^{(n-1)},\phi_1^{(n-1)},{\bf v_{1\perp}'}^{(n-1)}\right) = - J\left(\phi_1^{(n-1)},\omega_1^{(n-1)}\right) \dots\hspace{2cm}\nonumber\\
&\hspace{3cm}- 2 (\partial_y^2 - \partial_z^2) \left(\mathcal{R}\left(v_1'^{(n-1)} w_1'^{*(n-1)}\right)\right) - 2 \partial_y \partial_z \left(w_1'^{(n-1)} w_1'^{*(n-1)} - v_1'^{(n-1)} v_1'^{*(n-1)}\right),\\
&N_2\left({\bf v_{1\perp}'}^{(n-1)},u_0^{(n-1)},p_1'^{(n-1)}\right) = - i \alpha u_0^{(n-1)} {\bf v_{1\perp}'}^{(n-1)} - \nabla_{\perp} p_1'^{(n-1)},
\end{eqnarray}
and for which $p_1'^{(n-1)}$ and $\phi_1^{(n-1)}$ have already been evaluated in a preliminary step:
\begin{eqnarray}
&\phi_1^{(n-1)} = \left( \nabla_{\perp}^2 \right)^{-1} \omega_1^{(n-1)},\\
&p_1'^{(n-1)} = 2 i \alpha \left(\alpha^2 - \nabla_{\perp}^2\right)^{-1}  \left(v_1'^{(n-1)} \partial_y u_0^{(n-1)} + w_1'^{(n-1)} \partial_z u_0^{(n-1)}\right).
\end{eqnarray}
The convergence criterion $\mathcal{L}$ is obtained on the $L_2$-norm of the right hand side of equations (\ref{shearprec1})--(\ref{shearprec3}) multiplied by $(1 + \epsilon \, \triangle t_i)/\triangle t_i$ where $i=1$ for equations (\ref{shearprec1}), (\ref{shearprec2}) and $i=2$ for equation (\ref{shearprec3}).

We assume that the mean equations (\ref{shearred1}), (\ref{shearred2}) are of the same time as the doubly diffusive convection equations: their diffusion term is comparable in amplitude to the advective term.
We therefore assume that a good preconditioner for them is the mixed one and set $\triangle t_1 = Re$, such that the term $I - \epsilon \, \triangle t_1 \, L_{1j}$, with $j=1, 2$ in equations (\ref{shearprec1}), (\ref{shearprec2}) becomes $I - L_{1j}$.

%Simulations are run for different values of $\triangle t_2$ starting from the solution at $Re = 1000$.
%All these simulations eventually fail, sooner or later depending on the value of $\triangle t_2$.
%None of these simulations made it to $Re=2000$ and the highest Reynolds number reached is $Re \approx 1537.10$ obtained for $\triangle t_2 = 100$.
%Another batch of simulations was run from $Re = 1500$, in which all runs failed before $Re = 1950$.
%The efficiency of these simulations is sensitively different from the original runs from $Re = 1000$ and points towards the need for relatively small step sizes in the continuation method.

%% To enrich these observations, simulations were also run from $Re = 500$ and $Re = 2000$.
%% The results are shown in figure \ref{efficiency}.
%% %%%%%%%%%%%%%%%%%%%%%%%%%%
%% \begin{figure}
%% \begin{center}
%% \includegraphics[width=0.5\textwidth]{figlowerinit.eps}
%% \end{center}
%% \caption{Last continuation point computed successfully along the lower branch shown in figure \ref{initsol}.}
%% \label{efficiency}
%% \end{figure}
%% %%%%%%%%%%%%%%%%%%%%%%%%%%

\begin{figure}
\begin{center}
\includegraphics[width=\textwidth]{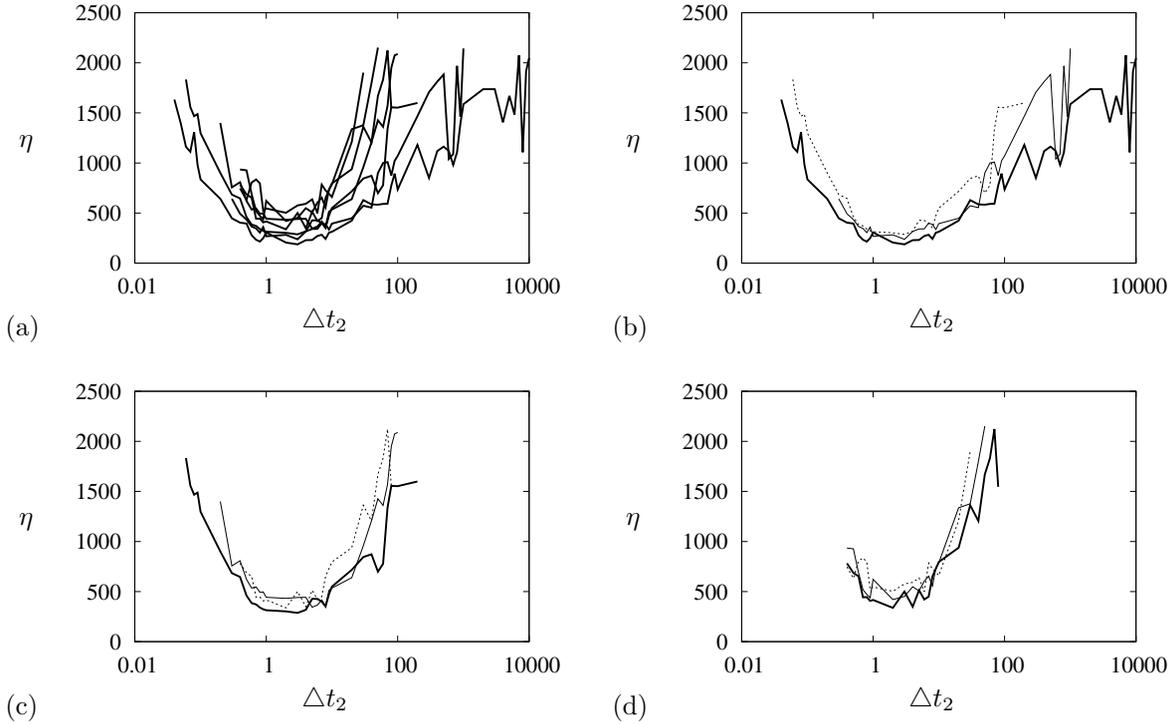}
\end{center}
\caption{(a) Number $\eta$ of gradient iterations needed to converge to a solution along the lower branch state (see figure \ref{initsol}) for $Re=500, 1000, 1500, 2000, 2500, 3000, 3500$ as a function of $\triangle t_2$ with $\triangle Re = 1$. (b) Subset of (a) for $Re=500$ (thick line), $Re=1000$ (thin line) and $Re=1500$ (dashed line). (c) Subset of (a) for $Re=1500$ (thick line), $Re=2000$ (thin line) and $Re=2500$ (dashed line). (d) Subset of (a) for $Re=2500$ (thick line), $Re=3000$ (thin line) and $Re=3500$ (dashed line).}
\label{wallylowerbranch1}
\end{figure}

To investigate the optimal preconditioner, we set up a number of simulations consisting in computing one continuation step with $\triangle Re = 1$.
These simulations are carried out for a range of $\triangle t_2$ with all other parameters unmodified.
The number of gradient iterations needed to converge is then recorded and reported in figure \ref{wallylowerbranch1} for some representative simulations.
%%%%%%%%%%%%%%%%%%%%%%%%%%
%%%%%%%%%%%%%%%%%%%%%%%%%%
As for the doubly diffusive problem, we observe that the algorithm fails at converging for small values of $\triangle t_2$.
This result is, in a sense, surprising: the fluctuation equations (\ref{shearred3}), (\ref{shearred4}) are weakly diffusive but still need preconditioning.
For large $\triangle t_2$, the method also fails, unlike for doubly diffusive convection.
More precisely, no simulation converged for $\triangle t_2 > 10^{5}$ and only a few successful events have been recorded for $\triangle t_2 > 500$ and $Re=500$.
Stokes preconditioning is hence not applicable here, as expected from the nature of the equations.

\begin{figure}
\begin{center}
\includegraphics[width=0.5\textwidth]{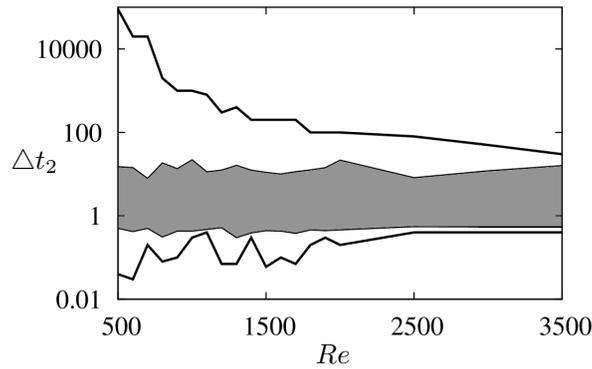}
\end{center}
\caption{Efficiency results for the continuation method in the ($\triangle t_2$,$Re$) plane. The outer lines indicate the working interval within which continuation works. The shaded region indicates the efficiency interval where continuation takes less than 200$\%$ of the number of gradient iterations needed at the optimal $\triangle t_2$ for the same value of $Re$.}
\label{figmarginslower}
\end{figure}

For relatively low values of $Re$, the continuation is rather permissive and a wide range of $\triangle t_2$ is allowed.
Increasing the Reynolds number has several effects on the continuation results.
First, the curves in figure \ref{wallylowerbranch1} move sensitively upwards, indicating that even for the optimal $\triangle t_2$, continuation requires more iterations from the conjugate gradient.
The boundaries between which the continuation method works also change dramatically.
The upper bound is $\triangle t_2 \approx 10^5$ at $Re=500$ and decreases abruptly until $\triangle t_2 \approx 100$ at $Re=2000$ before continuing to decrease, although less strongly.
The lower bound displays more irregularities, with values around $\triangle t_2 \approx 0.02$ for $Re=500$ and $Re=1500$.
Nonetheless, a trend is clearly observed from figure \ref{wallylowerbranch1}(a)--(d): the lower bound increases slowly until $\triangle t_2 \approx 0.4$ for $Re = 3500$.
As a result of the trend of both boundaries, the $\triangle t_2$ interval for which the continuation method works becomes narrower as $Re$ increases.
These tendencies are confirmed in figure \ref{figmarginslower}.
%%%%%%%%%%%%%%%%%%%%%%%%%%
%%%%%%%%%%%%%%%%%%%%%%%%%%

In addition to the parameter interval in which the continuation method works, it is important to know where it is efficient.
Indeed, a factor of $5$ can be observed between results at the optimum $\triangle t_2$ and those on the edge of the working interval (see figure \ref{wallylowerbranch1}).
Figure \ref{figmarginslower} reports the interval in which continuation is at least half as efficient as for the optimal $\triangle t_2$. 
This interval is calculated at each $Re$ in the following way: the optimal parameter $\triangle t_{2opt}$ for which the lowest number of gradient iterations is needed $\eta_{opt}$ is recorded and the interval is defined as the closest intersections from $\triangle t_{2opt}$ between the $\eta(\triangle t_2)$ curve and $\eta = 2 \eta_{opt}(Re)$.
Despite the fact that the working interval shrinks as $Re$ increases, the efficiency interval appears steady and spans approximately $0.5 \le \triangle t_2 \le 10$.
Note that the working interval shrinks principally from above with the Reynolds number due to the solution losing progressively its dissipative character.

\begin{figure}
\begin{center}
\includegraphics[width=\textwidth]{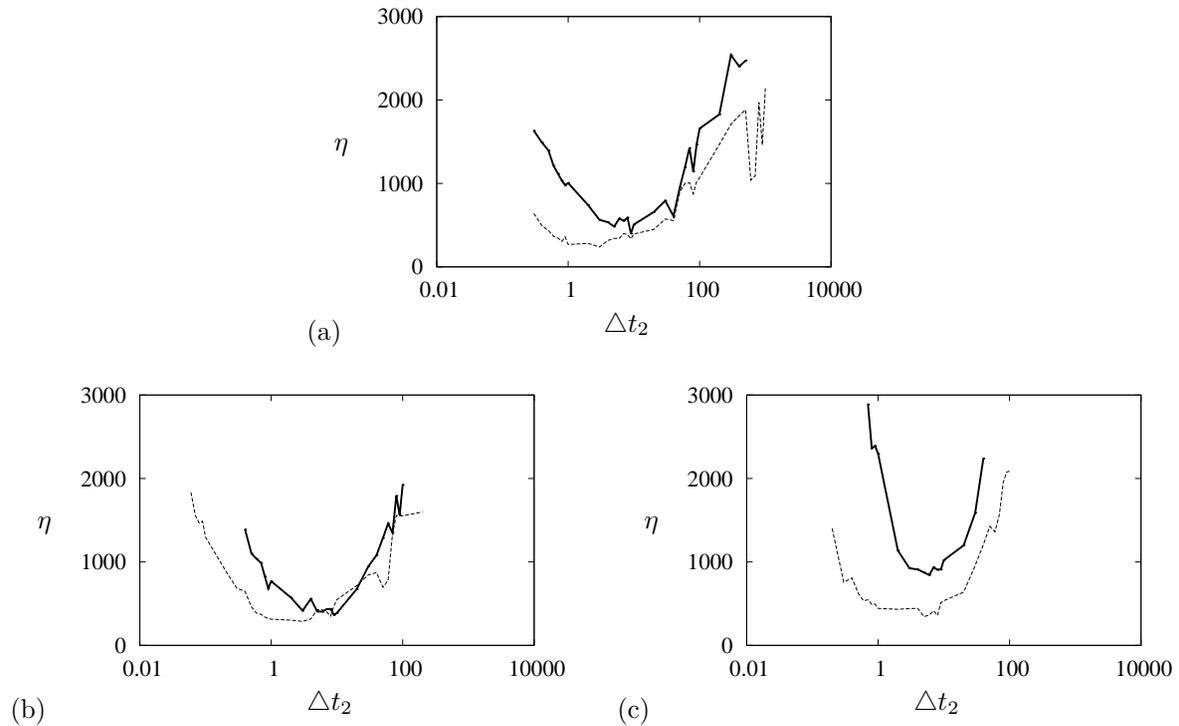}
\end{center}
\caption{Comparison of the continuation method efficiency between the lower branch solution (dashed lines, see also figure \ref{wallylowerbranch1}) and the upper branch solution (thick line). Results are shown for $Re = 1000$ (a), $Re=1500$ (b) and $Re=2000$ (c).}
\label{figupper}
\end{figure}

The above observations are compared to results along the upper branch at representative values of the Reynolds number in figure \ref{figupper}.
%%%%%%%%%%%%%%%%%%%%%%%%%%
%%%%%%%%%%%%%%%%%%%%%%%%%%
The continuation along the upper branch is substantially more complicated than along the lower branch.
For $Re = 1000$ and $Re = 1500$, one can see that for an optimal preconditioner at $\triangle t_{2opt}$, the upper branch is slightly more computationally demanding that the lower branch, however, away from the optimum, the upper branch solution is much harder to compute.
In addition to the number of gradient iterations $\eta$ that increases more abruptly for the upper branch away from $\eta_{opt}$, figure \ref{figupper} also indicates that the working interval of the preconditioner is narrower for the upper branch than for the lower branch.
These observations are confirmed and enhanced as the Reynolds number increases, as shown for $Re=2000$ in figure \ref{figupper}(c).

%%%%%%%%%%%%%%%%%%%%%%%%%%%%%%%%%%%%%%%%%%%%%%%%%%%%%%%%%%%%%%%%%%%%%%%
%%%%%%%%%%%%%%%%%%%%%%%%%%%%%%%%%%%%%%%%%%%%%%%%%%%%%%%%%%%%%%%%%%%%%%%
\section{Discussion}
%%%%%%%%%%%%%%%%%%%%%%%%%%%%%%%%%%%%%%%%%%%%%%%%%%%%%%%%%%%%%%%%%%%%%%%
%%%%%%%%%%%%%%%%%%%%%%%%%%%%%%%%%%%%%%%%%%%%%%%%%%%%%%%%%%%%%%%%%%%%%%%

In this paper, I have presented a preconditioner for numerical continuation of viscous incompressible stationary flows based on the Stokes preconditioner \cite{Mamun95}.
This preconditioner is easily constructed based on a time-stepper and has one parameter, $\triangle t$, corresponding to the time-step of the time-stepper.
In the small $\triangle t$ case, the preconditioner is weak and the limit $\triangle t \rightarrow 0$ corresponds to no preconditioning.
The case of large $\triangle t$ leads to Stokes preconditioning where the preconditioner approximates a Laplacian.
In the intermediate case, the preconditioner takes the form: $I-\triangle t L$, where $I$ is the identity, $\triangle t$ the parameter and $L$ the linear (Laplacian) operator.
The preconditioner is applied to two cases: three-dimensional doubly diffusive convection and a reduced model of shear flow.
In both cases, it is shown that the optimum preconditioner is mixed.

In the problem of doubly diffusive convection, the use of a preconditioner was necessary to carry out continuation successfully.
The Stokes preconditioner provides a good solution as it allows for continuation at a relatively constant and predictable speed but it is possible to improve it by using a mixed preconditioner.
The number of gradient iterations is then reduced by up to $60\%$ on the lower segment case ($140$ on average versus $350$ for Stokes preconditioning).

The reduced model of shear flow is more complex to deal with.
The set of equations comprising the reduced model is treated (simultaneously) using two different preconditioners.
A first subset of the model is highly diffusive and the mixed preconditioner $I - L$ has been successfully used.
The other subset requires preconditioning but unlike doubly diffusive convection, Stokes preconditioning does not work owing to the fact that the solutions studied are only weakly diffusive.
The optimum preconditioner is not only sensitive to the solution but also to the parameter $Re$: $I - \triangle t \, Re^{-1} \, L$ with $0.5 < \triangle t < 10$.

There exists only few numerical continuation methods in fluid systems \cite{Sanchez02,Viswanath09,channelflow} and the method presented here possesses some assets.
It is easy to implement as it is based on a first order impicit Euler time-scheme and it is adaptive through fine-tuning of the parameter $\triangle t$.
When searching for stationary solutions, the system of equations can be split into several sub-systems, each of which can then be treated with a different $\triangle t$ and parameter-dependent preconditioning is also easily implementable and proved efficient in one test case studied here.

\section*{Acknowledgments}
This author is grateful to L. S. Tuckerman for discussions and encouragements and thanks A. Bergeon for providing some of the routines used.


\begin{thebibliography}{10}

\bibitem{Allgower80}
{\sc E.~L. Allgower and K.~Georg}, {\em Simplicial and continuation methods for
  approximations, fixed points and solutions to systems of equations}, SIAM
  Rev., 22 (1980), pp.~28--85.

\bibitem{Allgower03}
\leavevmode\vrule height 2pt depth -1.6pt width 23pt, {\em Introduction to
  Numerical Continuation Methods}, Society for Industrial and Applied
  Mathematics, 2003.

\bibitem{Allgower85}
{\sc E.~L. Allgower and P.~H. Schmidt}, {\em An algorithm for piecewise-linear
  approximation of an implicitly defined manifold}, SIAM J. Numer. Anal., 22
  (1985), pp.~322--346.

\bibitem{Barkley02}
{\sc D.~Barkley, M.~G.~M. Gomes, and R.~D. Henderson}, {\em Three-dimensional
  instability in flow over a backward-facing step}, J. Fluid Mech., 473 (2002),
  pp.~167--190.

\bibitem{Batiste06}
{\sc O.~Batiste, E.~Knobloch, A.~Alonso, and I.~Mercader}, {\em Spatially
  localized binary-fluid convection}, J. Fluid Mech., 560 (2006), pp.~149--158.

\bibitem{Beaume13rc1}
{\sc C.~Beaume, A.~Bergeon, H.-C. Kao, and E.~Knobloch}, {\em Convectons in a
  rotating fluid layer}, J. Fluid Mech., 717 (2013), pp.~417--448.

\bibitem{Beaume11}
{\sc C.~Beaume, A.~Bergeon, and E.~Knobloch}, {\em Homoclinic snaking of
  localized states in doubly diffusive convection}, Phys. Fluids, 23 (2011),
  p.~094102.

\bibitem{Beaume13ddd1}
\leavevmode\vrule height 2pt depth -1.6pt width 23pt, {\em Convectons and
  secondary snaking in three-dimensional natural doubly diffusive convection},
  Phys. Fluids, 25 (2013), p.~024105.

\bibitem{Beaume15pre}
{\sc C.~Beaume, G.~P. Chini, K.~Julien, and E.~Knobloch}, {\em Reduced
  description of exact coherent states in parallel shear flows}, Phys. Rev. E,
  91 (2015), p.~043010.

\bibitem{Beaume13rc2}
{\sc C.~Beaume, H.-C. Kao, E.~Knobloch, and A.~Bergeon}, {\em Localized
  rotating convection with no-slip boundary conditions}, Phys. Fluids, 25
  (2013), p.~124105.

\bibitem{Beaume13ddd2}
{\sc C.~Beaume, E.~Knobloch, and A.~Bergeon}, {\em Nonsnaking doubly diffusive
  convectons and the twist instability}, Phys. Fluids, 25 (2013), p.~114102.

\bibitem{Bergeon98}
{\sc A.~Bergeon, D.~Henry, H.~{Ben Hadid}, and L.~S. Tuckerman}, {\em Marangoni
  convection in binary mixtures with soret effect}, J. Fluid Mech., 375 (1998),
  pp.~143--177.

\bibitem{Bergeon02}
{\sc A.~Bergeon and E.~Knobloch}, {\em Natural doubly diffusive convection in
  three-dimensional enclosures}, Phys. Fluids, 14 (2002), pp.~3233--3250.

\bibitem{Boronska10}
{\sc K.~Boro{\'n}ska and L.~S. Tuckerman}, {\em Extreme multiplicity in
  cylindrical {R}ayleigh--{B}\'enard convection. {II}. {B}ifurcation diagram
  and symmetry classification}, Phys. Rev. E, 81 (2010), p.~036321.

\bibitem{Chantry14}
{\sc M.~Chantry, A.~P. Willis, and R.~R. Kerswell}, {\em Genesis of
  streamwise-localized solutions from globally periodic traveling waves in pipe
  flow}, Phys. Rev. Lett., 112 (2014), p.~164501.

\bibitem{Clever97}
{\sc R.~M. Clever and F.~H. Busse}, {\em Tertiary and quaternary solutions for
  plane couette flow}, J. Fluid Mech., 344 (1997), pp.~137--153.

\bibitem{Clewley07}
{\sc R.~H. Clewley, W.~E. Sherwood, M.~D. La{M}ar, and J.~M. Guckenheimer},
  {\em Py{DST}ool, a software environment for dynamical systems modeling},
  (2007).

\bibitem{Vandervorst92}
{\sc H.~A.~Van der Vorst}, {\em Bi-cgstab: A fast and smoothly converging
  variant of bi-cg for the solution of nonsymmetric linear systems}, SIAM J.
  Sci. and Stat. Comput., 13 (1992), pp.~631--644.

\bibitem{Dijkstra14}
{\sc H.~A. Dijkstra, F.~W. Wubs, A.~K. Cliffe, E.~Doedel, I.~F. Dragomirescu,
  B.~Eckhardt, A.~Y. Gelfgat, A.~L. Hazel, V.~Lucarini, A.~G Salinger, E.~T.
  Phipps, J.~Sanchez-Umbria, H.~Schuttelaars, L.~S. Tuckerman, and U.~Thiele},
  {\em Numerical bifurcation methods and their applications to fluid dynamics:
  analysis beyond simulation}, Commun. Comput. Phys., 15 (2014), pp.~1--45.

\bibitem{Doedel08}
{\sc E.~J. Doedel, A.~R. Champneys, F.~Dercole, T.~Fairgrieve, Y.~Kuznetsov,
  B.~Oldeman, R.~Paffenroth, B.~Sandstede, X.~Wang, and C.~Zhang}, {\em
  {AUTO-07P: Continuation and Bifurcation Software for Ordinary Differential
  Equations}}, Feb. 2008.

\bibitem{Engelborghs01}
{\sc K.~Engelborghs, T.~Luzyanina, and G.~Samaey}, {\em {DDE-BIFTOOL} v. 2.00:
  a {M}atlab package for bifurcation analysis of delay differential equations},
  Technical Report TW-330, Dep. Comp. Sci., KU Leuven,  (2001).

\bibitem{Faisst03}
{\sc H.~Faisst and B.~Eckhardt}, {\em Traveling waves in pipe flow}, Phys. Rev.
  Lett., 91 (2003), p.~224502.

\bibitem{Feudel11}
{\sc F.~Feudel, K.~Bergemann, L.~S. Tuckerman, C.~Egbers, B.~Futterer,
  M.~Gellert, and R.~Hollerbach}, {\em Convection patterns in a spherical fluid
  shell}, Phys. Rev. E, 83 (2011), p.~046304.

\bibitem{Frigo05}
{\sc M.~Frigo and S.~G. Johnson}, {\em The design and implementation of
  {FFTW3}}, Proceedings of the IEEE, 93 (2005), pp.~216--231.

\bibitem{channelflow}
{\sc J.~F. Gibson}, {\em {Channelflow}: {A} spectral {Navier-Stokes} simulator
  in {C}++}, tech. report, U. New Hampshire, 2012.
\newblock {\tt {Channelflow.org}}.

\bibitem{Gibson14}
{\sc J.~F. Gibson and E.~Brand}, {\em Spanwise-localized solutions of planar
  shear flows}, J. Fluid Mech., 745 (2014), pp.~25--61.

\bibitem{Gibson09}
{\sc J.~F. Gibson, J.~Halcrow, and P.~Cvitanovi{\'c}}, {\em Equilibrium and
  travelling-wave solutions of plane couette flow}, J. Fluid Mech., 638 (2009),
  pp.~243--266.

\bibitem{Henry00}
{\sc D.~B. Henry and A.~Bergeon}, {\em Continuation methods in fluid dynamics},
  Vieweg, 2000.

\bibitem{Karniadakis91}
{\sc G.~E. Karniadakis, M.~Israeli, and S.~A. Orszag}, {\em High-order
  splitting methods for the incompressible navier--stokes equations}, J. Comp.
  Phys., 97 (1991), pp.~414--443.

\bibitem{Kawahara12}
{\sc G.~Kawahara, M.~Uhlmann, and L.~van Veen}, {\em The significance of simple
  invariant solutions in turbulent flows}, Annu. Rev. Fluid Mech., 44 (2012),
  pp.~203--225.

\bibitem{Keller77}
{\sc H.-B. Keller}, {\em Numerical solutions of bifurcation and non-linear
  eigenvalues problem: {A}pplication of bifurcation theory}, Academic Press New
  York,  (1977).

\bibitem{Krauskopf07}
{\sc B.~Krauskopf, H.~M. Osinga, and J.~Gal{\'a}n-{V}iosque}, {\em Numerical
  continuation methods for dynamical systems: {P}ath following and boundary
  value problems}, Springer-Verlag, 2007.

\bibitem{Kuznetsov03}
{\sc Y.~A. Kuznetsov, A.~Dhooge, and W.~Govaerts}, {\em Matcont: a {MATLAB}
  package for numerical bifurcation analysis of {ODE}s}, ACM Trans. Math.
  Softw., 29 (2003), pp.~141--164.

\bibitem{Lojacono11}
{\sc D.~{Lo Jacono}, A.~Bergeon, and E.~Knobloch}, {\em Magnetohydrodynamic
  convectons}, J. Fluid Mech., 687 (2011), pp.~595--605.

\bibitem{Lojacono13}
\leavevmode\vrule height 2pt depth -1.6pt width 23pt, {\em Three-dimensional
  spatially localized binary-fluid convection in a porous medium}, J. Fluid
  Mech., 730 (2013), p.~R2.

\bibitem{Mamun95}
{\sc C.~K. Mamun and L.~S. Tuckerman}, {\em Asymmetry and hopf bifurcation in
  spherical couette flow}, Phys. Fluids, 7 (1995), pp.~80--91.

\bibitem{Melnikov14}
{\sc K.~Melnikov, T.~Kreilos, and B.~Eckhardt}, {\em Long-wavelength
  instability of coherent structures in plane couette flow}, Phys. Rev. E, 89
  (2014), p.~043008.

\bibitem{Mercader09}
{\sc I.~Mercader, O.~Batiste, A.~Alonso, and E.~Knobloch}, {\em Localized
  pinning states in closed containers: {H}omoclinic snaking without
  bistability}, Phys. Rev. E, 80 (2009), p.~025201(R).

\bibitem{Mercader11}
\leavevmode\vrule height 2pt depth -1.6pt width 23pt, {\em Convectons,
  anticonvectons and multiconvectons in binary fluid convection}, J. Fluid
  Mech., 667 (2011), pp.~586--606.

\bibitem{Nagata90}
{\sc M.~Nagata}, {\em Three-dimensional finite-amplitude solutions in plane
  {C}ouette flow: bifurcation from infinity}, J. Fluid Mech., 217 (1990),
  pp.~519--527.

\bibitem{Pringle07}
{\sc C.~C.~T. Pringle and R.~R. Kerswell}, {\em Asymmetric, helical and
  mirror-symmetric traveling waves in pipe flow}, Phys. Rev. Lett., 99 (2007),
  p.~074502.

\bibitem{Saad03}
{\sc Yousef Saad}, {\em {Iterative Methods for Sparse Linear Systems, Second
  Edition}}, Society for Industrial and Applied Mathematics, 2003.

\bibitem{Sanchez02}
{\sc J.~Sanchez, F.~Marques, and J.~M. Lopez}, {\em A continuation and
  bifurcation technique for navier--{S}tokes flows}, J. Comput. Phys., 180
  (2002), pp.~78--98.

\bibitem{Schneider10}
{\sc T.~M. Schneider, J.~F. Gibson, and J.~Burke}, {\em Snakes and ladders:
  {L}ocalized solutions of plane couette flow}, Phys. Rev. Lett., 104 (2010),
  p.~104501.

\bibitem{Schneider08}
{\sc T.~M. Schneider, J.~F. Gibson, M.~Lagha, F.~De Lillo, and B.~Eckhardt},
  {\em Laminar-turbulent boundary in plane couette flow}, Phys. Rev. E, 78
  (2008), p.~037301.

\bibitem{Seydel09}
{\sc R.~Seydel}, {\em Practical Bifurcation and Stability Analysis},
  Interdisciplinary Applied Mathematics, Springer, 2009.

\bibitem{Seydel87}
{\sc R.~Seydel and V.~Hlavacek}, {\em Role of continuation in engineering
  analysis}, Chem. Eng. Sci., 42 (1987), pp.~1281--1295.

\bibitem{Torres14}
{\sc J.~F. Torres, D.~Henry, A.~Komiya, and S.~Maruyama}, {\em Bifurcation
  analysis of steady natural convection in a tilted cubical cavity with
  adiabatic sidewalls}, J. Fluid Mech., 756 (2014), pp.~650--688.

\bibitem{Torres13}
{\sc J.~F. Torres, D.~Henry, A.~Komiya, S.~Maruyama, and H.~{Ben Hadid}}, {\em
  Three-dimensional continuation study of convection in a tilted rectangular
  enclosure}, Phys. Rev. E, 88 (2013), p.~043015.

\bibitem{Tuckerman89}
{\sc L.~S. Tuckerman}, {\em Steady-state solving via stokes preconditioning;
  recursion relations for elliptic operators}, in 11th International Conference
  on Numerical Methods in Fluid Dynamics, D.L. Dwoyer, M.Y. Hussaini, and R.G.
  Voigt, eds., vol.~323 of Lecture Notes in Physics, Springer Berlin
  Heidelberg, 1989, pp.~573--577.

\bibitem{Uecker14}
{\sc H.~Uecker, D.~Wetzel, and J.~D.~M. Rademacher}, {\em Pde2path -- a
  {M}atlab package for continuation and bifurcation in 2{D} elliptic systems},
  Numerical Mathematics: Theory, Methods and Applications, 7 (2014),
  pp.~58--106.

\bibitem{Viswanath09}
{\sc D.~Viswanath}, {\em The critical layer in pipe flow at high reynolds
  numbers}, Phil. Trans. R. Soc. A, 367 (2009), pp.~561--576.

\bibitem{Waleffe03}
{\sc F.~Waleffe}, {\em Homotopy of exact coherent structures in plane shear
  flows}, Phys. Fluids, 15 (2003), pp.~1517--1534.

\bibitem{Wang07}
{\sc J.~Wang, J.~Gibson, and F.~Waleffe}, {\em Lower branch coherent states in
  shear flows: {T}ransition and control}, Phys. Rev. Lett., 98 (2007),
  p.~204501.

\bibitem{Wedin04}
{\sc H.~Wedin and R.~R. Kerswell}, {\em Exact coherent structures in pipe flow:
  travelling wave solutions}, J. Fluid Mech., 508 (2004), pp.~333--371.

\end{thebibliography}
\end{document}